\documentstyle[amssymb,12pt]{amsart}

\ifx\mathrm\undefined\let\mathrm\bf\fi
\ifx\mathbf\undefined\let\mathbf\bf\fi

 \def\Sum{\sum\limits}

\let\epe\epsilon \let\eps\varepsilon \let\epsilon\eps

\let\la\lambda \let\La\Lambda
 
 \let\phi\varphi

\newcommand{\I}{{{\frak I}\,}}
\newcommand{\G}{{{\frak G}\,}}
\newcommand{\F}{{{\frak F}\,}}
\newcommand{\Zb}{{{\mathcal Z}}}

\newcommand{\Z}{{\Bbb Z}}

\newcommand{\R}{{\Bbb R}}  
\newcommand{\C}{{\Bbb C}}  
  
\newcommand{\Ref}[1]{{$($\ref{#1}$)$}}
\newcommand{\bean}{\begin{eqnarray}}
\newcommand{\eean}{\end{eqnarray}}
\newcommand{\be}{\begin{displaymath}}
\newcommand{\ee}{\end{displaymath}}
\newcommand{\bea}{\begin{eqnarray*}}
\newcommand{\eea}{\end{eqnarray*}}

\newcommand{\res}{{\operatorname{res}}}

\newcommand{\T}{\otimes}

\newcommand{\vs}{\vspace{1.5\baselineskip}}

\newenvironment{proof}{\noindent{\it Proof\/}:\rm}{$\;\Box$\par\vs}

\newtheorem
{thm}{Theorem}

\newtheorem
{lemma}[thm]{Lemma}

\newtheorem
{corollary}[thm]{Corollary}

\newcommand{\End}{{\operatorname{End}}}

\newcommand{\Imag}{{\operatorname{Im}}\,}
\newcommand{\Real}{{\operatorname{Re}}\,}

\begin{document}

\title[The qKZ Equation in Tensor Products of Irreducible Modules]
{ The Quantized Knizhnik-Zamolodchikov Equation in Tensor Products of Irreducible $sl_2$-Modules}
\author{E. Mukhin and A. Varchenko}
\maketitle
\vskip-.5\baselineskip
\centerline{${}^*${\it Department of Mathematics,
University of North Carolina at Chapel Hill,}}
\centerline{\it Chapel Hill, NC 27599-3250, USA}
\centerline{{\it E-mail addresses:} {\rm mukhin@@math.unc.edu,
av@@math.unc.edu}}
\bigskip
\medskip
\centerline{August, 1997}
\bigskip
\medskip

\begin{abstract}

We consider the quantized Knizhnik-Zamolodchikov
difference equation (qKZ)  with values in a tensor product of irreducible
$sl_2$ modules, the equation defined in terms of rational R-matrices.
We solve the equation in terms of multidimensional q-hypergeometric
integrals. We identify  the space of solutions of the qKZ equation
with the tensor product of the corresponding modules over the quantum
group $U_qsl_2$. We compute the monodromy of the qKZ equation in terms
of the trigonometric $R$-matrices.

\end{abstract}

\section{Introduction.}
\subsection{}
In this paper we solve the rational quantized Knizhnik-Zamolodchikov difference
equation ( qKZ ) with values in a tensor product of irreducible highest weight
$sl_2$-modules. The rational qKZ equation 
is a system of difference equations for a function $\Psi(z_1,...,z_n)$ with 
values in a tensor product $M_1\otimes ...\otimes M_n$ of $sl_2$ modules.
The system of equations has the form
\bea
\Psi(z_1,\dots,z_m+p,\dots,z_n)=R_{m, {m-1}}(z_m-z_{m-1}+p)\dots
R_{m,1}(z_m-z_1+p)e^{-\mu h_m}\times
\\
\times R_{m,n}(z_m-z_n)\dots R_{m,{m+1}}(z_m-z_{m+1})\Psi(z_1,\dots,z_n),
\eea
 $m=1,\dots,n$, where $p$, $\mu$ are complex parameters of the qKZ equation, 
$h$ is a generator of the Cartan subalgebra of $sl_2$, 
$h_m$ is the operator $h$ acting in the $m$-th factor,
$R_{i,j}(x)$ is the rational $R$-matrix 
$R_{M_iM_j}(x)\in \End(M_i\T M_j)$ acting in the $i$-th and $j$-th factors
of the tensor product. In this paper we consider only steps $p$ with negative real part.

The qKZ equation is an important system of 
difference equations. The qKZ equation was introduced in \cite{FR} as an equation
for matrix elements of vertex operators of a quantum affine algebra.
An important special case of the qKZ equation had been introduced earlier in \cite{S}
as equations for form factors in integrable quantum field theory.
Later, the qKZ equation was derived as an equation for correlation functions
in lattice integrable models, cf. \cite{JM} and references therein.

Solutions of the rational qKZ equation with values in a tensor product of
$sl_2$ Verma modules $V_{\la_1}\T ...\T V_{\la_n}$ with generic highest weights 
$\la_1,\dots,\la_n$
were constructed in \cite{TV1}. The solutions have the form
\bea
\Psi(z)\,=\, \sum_{k_1,...,k_n} I_{k_1,...,k_n}(z) f^{k_1}v_1 \otimes
... \otimes  f^{k_n}v_n,
\eea
where $\{f^{k_1}v_1\T\ldots\T f^{k_n}v_n\}\in V_{\la_1}\T ...\T V_{\la_n}$ is 
the standard basis in the tensor product of $sl_2$ Verma modules,
and the coefficients $I_{k_1,...,k_n}(z)$ are given by suitable multidimensional
q-hypergeometric integrals. 

The space of solutions
of the qKZ equation with values in a tensor product of
$sl_2$ Verma modules with generic highest weights 
was described in \cite{TV1} in terms of the representation theory of the quantum group
$U_q(sl_2)$ with $q=e^{\pi i/p}$. Namely, consider the tensor product 
$V^q_{\la_1}\T ...\T V^q_{\la_n}$ of $U_q(sl_2)$ Verma modules,
where $V_{\la_j}^q$ is the deformation of the $sl_2$ Verma module $V_{\la_j}$.
It was shown in \cite{TV1} that there is a natural isomorphism of the space $\cal S$
of meromorphic solutions of the qKZ equation
with values in the tensor product
$V_{\la_1}\T ...\T V_{\la_n}$ and the space
$V^q_{\la_1}\T ...\T V^q_{\la_n}\T F$, where $F$ is the space of
meromorphic functions in $z_1,...,z_n$, $p$-periodic with respect to each of the variables,
\be
 V^q_{\la_1}\T ...\T V^q_{\la_n}\T F \,\simeq \, \cal S.
\ee
This isomorphism was used in \cite{TV1} to compute asymptotic solutions of the qKZ equation
with values in a tensor product of
$sl_2$ Verma modules with generic highest weights and to compute the
transition functions
between asymptotic solutions
in terms of the trigonometric $R$-matrices acting in $V^q_{\la_1}\T ...\T V^q_{\la_n}$.

Assume that the Verma module $V_{\la_j}$ is reducible and
vectors $\{f^kv_j\}_{k\ge N_j}$ generate a proper submodule $S_{\la_j}\subset V_{\la_j}$.
Then the vectors $\{f^kv_j\}_{k< N_j}$ form a basis in the irreducible
module $L_{\la_j}=V_{\la_j}/S_{\la_j}$. If $V_{\la_j}$ is irreducible, then we set
$N_j=\infty$.

Return to the sum $\Psi(z)$.
If the Verma modules of the tensor product $V_{\la_1}\T ...\T V_{\la_n}$
become reducible, then some of the coefficients in the sum become divergent.
In this paper we show that the restricted sum
\bea
\Psi^0(z)\,=\, 
\sum_{k_1<N_1 ,...,k_n<N_n} I_{k_1,...,k_n}(z) 
f^{k_1}v_1 \T\ldots\T f^{k_n}v_n
\eea
remains well defined even when some of the Verma modules become reducible. 

Moreover, we show that the sum $\Psi^0(z)$
defines a solution of the qKZ equation with values
in the tensor product $L_{\la_1}\T ... \T L_{\la_n}$ of
irreducible $sl_2$ modules, and under certain conditions all solutions 
have this form.

These results allow us to describe the space of solutions to the qKZ equation with values
in $L_{\la_1}\T\ldots\T L_{\la_n}$ 
in terms of representation theory of the quantum group
$U_q(sl_2)$.
Namely, consider the tensor product
$L^q_{\la_1}\T ...\T L^q_{\la_n}$ of $U_qsl_2$ modules,
where $L_{\la_j}^q$ is the deformation of the $sl_2$ irreducible module 
$L_{\la_j}$. We show that there is a natural isomorphism of the space ${\cal S}$
of meromorphic solutions of the qKZ equation
with values in the tensor product
$L_{\la_1}\T ...\T L_{\la_n}$ and the space
$L^q_{\la_1}\T ...\T L^q_{\la_n}\T F$, 
\be  
 L^q_{\la_1}\T ...\T L^q_{\la_n}\T F \,\simeq \, {\cal S}.
\ee

We compute asymptotic solutions of the qKZ equation
with values in a tensor product of irreducible
$sl_2$ Verma modules and the
transition functions between the asymptotic solutions. The transition functions are given
in terms of the trigonometric $R$-matrices acting in $L^q_{\la_1}\T ...\T L^q_{\la_n}$.

In this paper we consider the rational qKZ equation associated with $sl_2$.
There are other types of the qKZ equation: the trigonometric qKZ equation \cite{FR}, 
\cite{TV2} and the elliptic qKZB equation \cite{F1},\cite{F2}, \cite{FTV1}, \cite{FTV2}.  
The trigonometric qKZ equation with values in a tensor product
of $U_q(sl_2)$ Verma modules with generic highest weights and the elliptic qKZB equation
with values in a tensor product
of $E_{\tau, \eta}(sl_2)$ Verma modules with generic highest weights were solved
in \cite{TV2} and \cite{FTV2}, respectively. Here $E_{\tau,\eta}(sl_2)$ is the elliptic 
quantum 
group associated to $sl_2$. In the next paper we will extend our results
to the trigonometric qKZ and elliptic qKZB equations
and describe solutions of these equations
with values in irreducible finite dimensional modules.

The paper is organized as follows.

In Section~\ref{general} we recall some facts about $sl_2$, $U_q(sl_2)$ 
and their representations. We define the rational $R$-matrix and 
the qKZ equation. In Section~\ref{pairing}
we describe integral representations of solutions of the qKZ equation
with values in a tensor product of $sl_2$ Verma modules.
The statements of results are given in Section~\ref{main}.
The proofs are collected in  Section~\ref{proofs}.

The authors thank V.Tarasov for useful discussions.

\section{General definitions and notations}\label{general}

\subsection{The Lie algebra $sl_2$}

Let $e$, $f$, $h$ be generators of the Lie algebra $sl_2$ such that 
\be
[h,e]=e, \qquad [h,f]=-f, \qquad [e,f]=2h.
\ee
For an $sl_2$ module $M$, 
let $M^*$ be its restricted dual with an $sl_2$ module structure defined by
\be
\langle e\varphi,x\rangle=\langle \varphi,fx\rangle,  \qquad
\langle f\varphi,x\rangle=\langle \varphi,ex\rangle,  \qquad
\langle h\varphi,x\rangle=\langle \varphi,hx\rangle
\ee
for all $x\in M$, $\varphi \in M^*$. The module $M^*$ is called the \emph{dual} 
module. 

For $\la \in \C$, 
denote $V_\la$ the $sl_2$  Verma module with  highest weight $\la$. Then
$V_\la=\bigoplus\limits_{i=0}^\infty \C f^iv$, where $v$ is a highest weight
vector. Denote $L_\la$ the irreducible module with highest weight $\la$.

Let $\La^+=\{0,\frac{1}{2},1,\frac{3}{2},2,...\}$ be the set of dominant 
weights. If $\la\in\La^+$, then $L_\la$ is  a $(2\la+1)$-dimensional module and
\be
L_\la\,\simeq\, V_\la/S_\la,
\ee
where $S_\la =\bigoplus\limits_{i=2\la+1}^\infty \C f^iv\subset
V_\la$ is the 
maximal proper submodule.  The
vectors $f^iv$, $i=0,\dots,2\la$, generate a
basis in $L_\la$.

For $\la\not\in\La^+$, $L_\la =V_\la$.
It is convenient to introduce $S_\la$ to be the zero submodule of $V_\la$, then
$L_\la\simeq V_\la/S_\la$, as we have for $\la\in\La^+$.

For an $sl_2$ module $M$ with  highest weght $\la$, denote by $(M)_l$ the 
subspace of weight $\la-l$, by $(M)^{sing}$ the kernel of the operator $e$, and 
by $(M)_l^{sing}$ the subspace $(M)_l\bigcap (M)^{sing}$.

\subsection{The algebra $U_qsl_2$.} 

Let q be a complex number different from $\pm 1$. 
Let $e_q$, $f_q$, $q^h$, $q^{-h}$ be generators of $U_qsl_2$ such that
\be
q^hq^{-h}=q^{-h}q^h=1, \qquad [e_q,f_q]=\frac{q^{2h}-q^{-2h}}{q-q^{-1}},
\ee
\be
q^he_q=qe_qq^h, \qquad  q^hf_q=q^{-1}f_qq^h.
\ee

A comultiplication $\Delta : U_qsl_2\to U_qsl_2\T U_qsl_2$ is given by
\be
\Delta(q^h)=q^h\T q^h,\qquad \Delta(q^{-h})=q^{-h}\T q^{-h},
\ee
\be
\Delta(e_q)=e_q\T q^h + q^{-h}\T e_q,\qquad 
\Delta(f_q)=f_q\T q^h + q^{-h}\T f_q.
\ee
The comultiplication defines a module structure on tensor products of 
$U_qsl_2$ modules.

For $\la \in \C$, denote $V_\la^q$ the $U_qsl_2$ Verma module with highest weight 
$q^\la$. Then
$V^q_\la=\bigoplus\limits_{i=0}^\infty \C f_q^iv^q$, where $v^q$ is a highest weight
vector. 

For $\la\in \La^+$, $S^q_\la=\bigoplus\limits_{i=2\la+1}^\infty \C f_q^iv^q$ is a submodule in $V^q_\la$.
Denote $L^q_\la$ the quotient module $V^q_\la/S^q_\la$. The module
$L^q_\la$ is the $(2\la+1)$-dimensional  
highest weight  module with highest weight $q^\la$. The 
vectors $f_q^iv^q$, $i=0,\dots,2\la$, generate a 
basis in $L^q_\la$.

For $\la\not\in\La^+$, let $L^q_\la=V^q_\la$.
It is convenient to introduce $S^q_\la$ to be the zero submodule of $V^q_\la$, then
$L^q_\la\simeq V^q_\la/S^q_\la$, as we have for $\la\in\La^+$.

For an $U_qsl_2$ module $M^q$ with highest weight $q^\la$, denote by 
$(M^q)_l$ the
subspace of weight $q^{\la-l}$, by $(M^q)^{sing}$ the kernel of the 
operator $e_q$, and
by $(M^q)_l^{sing}$ the subspace $(M^q)_l\bigcap (M^q)^{sing}$.

\subsection{The rational $R$-matrix}\label{rat $R$-matrix} 

The \emph{Yangian} $Y(sl_2)$ is a Hopf algebra which has a family of homomorphisms to the universal
enveloping algebra of $sl_2$, 
$Y(sl_2)\to U(sl_2)$, depending on a complex parameter. 
Therefore, each $sl_2$ module $M$ carries a Yangian module structure $M(x)$ depending on a parameter,
see \cite{CP},\cite{TV1}.

For $sl_2$ irreducible highest weight modules $L_{\la_1},L_{\la_2}$ and generic numbers $x,y\in\C$, the
Yangian modules $L_{\la_1}(x)\T L_{\la_2}(y)$
and  $L_{\la_2}(y)\T L_{\la_1}(x)$ are irreducible and isomorphic. 
There exists a unique intertwiner which sends $v_1\T v_2$ to $v_2\T v_1$, where $v_i$ is a highest
vector
in $L_{\la_i}$, $i=1,2$. This intertwiner has the form $PR_{L_{\la_1}L_{\la_2}}(x-y)$, where $P$ is the
operator of permutation of factors. The operator
$R_{L_{\la_1}L_{\la_2}}(x)\in \End(L_{\la_1}\T L_{\la_2})$ is called the
\emph{rational $R$-matrix}, see details in \cite{TV1}. The rational
$R$-matrix commutes
with the $sl_2$ action on the tensor product $L_{\la_1}\T L_{\la_2}$.

Let $\la_1,\la_2\not\in \La^+$. Let 
\be
V_{\la_1}\T
V_{\la_2}=\bigoplus\limits_{l=0}^{\infty} V_{\la_1+\la_2-l}
\ee
be the decomposition of the tensor product of $sl_2$ Verma modules into the 
direct sum of irreducibles, and let
$\Pi^{(l)}$ be the projector on $V_{\la_1+\la_2-l}$ along the other summands.
  
There is a formula for the rational $R$-matrix, 
\be
R_{V_{\la_1}V_{\la_2}}(x)=\sum_{l=0}^\infty\Pi^{(l)}\prod_{s=0}^{l-1}
\frac{x+\la_1+\la_2-s}{x-\la_1-\la_2+s},
\ee
see \cite{KRS},\cite{T}.

Let $\la_1,\la_2\in \La^+$. Let
\be
L_{\la_1}\T L_{\la_2}=\bigoplus\limits_{l=0}^{2\min\{\la_1,\la_2\}}L_{\la_1+\la_2-l}
\ee
be the decomposition of the tensor product of finite dimensional irreducible $sl_2$ modules into
the direct sum of irreducibles, and let
$\Pi^{(l)}$ be the projector on $L_{\la_1+\la_2-l}$ along the other summands.

There is a formula for the rational $R$-matrix,
\be
R_{L_{\la_1}L_{\la_2}}(x)=\sum_{l=0}^{2\min\{\la_1,\la_2\}}\Pi^{(l)}\prod_{s=0}^{l-1}
\frac{x+\la_1+\la_2-s}{x-\la_1-\la_2+s},
\ee
see \cite{KRS},\cite{T}.

The vector spaces $V_{\la_1}\T V_{\la_2}$ for different values of $\la_1,\la_2$ are identified by
distinguished bases $\{f^{l_1}v_1\T f^{l_2}v_2\,|\, l_1,l_2\in\Z_{\ge 0}\}$.

\begin{thm}\label{rational $R$-matrix}
1. The rational $R$-matrix $R_{V_{\la_1}V_{\la_2}}(x)\in \End(V\T V)$ is a meromorphic
function of $x,\la_1,\la_2$. The poles of $R_{V_{\la_1}V_{\la_2}}(x)$
 have the form 
$x-\la_1-\la_2+s=0$, where $s\in\Z_{\ge 0}$.

2. Let $x$ be generic. Then the rational $R$-matrix $R_{V_{\la_1}V_{\la_2}}(x)$ preserves 
$S_{\la_1}\T V_{\la_2}+V_{\la_1}\T S_{\la_2}$.

3. Let $V_{\la_1}\T V_{\la_2}\to L_{\la_1}\T L_{\la_2}$ be the canonical factorization map. Then
for
generic $x$, the rational $R$-matrix  $R_{V_{\la_1}V_{\la_2}}(x)$ can be factorized to an operator
$R(x)\in \End(L_{\la_1}\T L_{\la_2})$ and, moreover,
$R(x)=R_{L_{\la_1}L_{\la_2}}(x)$. 
\end{thm}
An analog of Theorem~\ref{rational $R$-matrix} for the elliptic $R$-matrix is proved in Theorems 8, 31 in \cite{FTV1}.
The same proof works for the rational $R$-matrix. We give an alternative proof in Section~\ref{proofs}.

\subsection{The trigonometric $R$-matrix}\label{trig $R$-matrix}

Let $q$ be a complex 
number, and not a root of unity.
The \emph{quantum affine algebra} $\widehat{U_qsl_2}$ is a Hopf algebra
which has a family of homomorphisms,
$\widehat{U_qsl_2}\to U_qsl_2$, depending on a complex parameter.
Therefore, each $U_qsl_2$ module $M^q$ carries a $\widehat{U_qsl_2}$ module structure $M^q(x)$
depending
on
a parameter, see  \cite{T}, \cite{CP}.

For $U_qsl_2$ irreducible highest weight modules $L^q_{\la_1},L^q_{\la_2}$ and generic numbers $x,y\in\C$, the
$\widehat{U_qsl_2}$ modules $L^q_{\la_1}(x)\T L^q_{\la_2}(y)$
and  $L^q_{\la_2}(y)\T L^q_{\la_1}(x)$ are isomorphic.
There exists a unique intertwiner which sends $v^q_1\T v^q_2$ to $v^q_2\T v^q_1$, where $v^q_i$ is a
highest vector
in $L^q_{\la_i}$, $i=1,2$. This intertwiner has the form $PR^q_{L^q_{\la_1}L^q_{\la_2}}(x/y)$, where $P$
is the
operator of permutation of factors. The operator $R^q_{L^q_{\la_1}L^q_{\la_2}}(x)\in
\End(L^q_{\la_1}\T L^q_{\la_2})$ is called the
\emph{trigonometric $R$-matrix}, see details in \cite{TV1}. 
The trigonometric $R$-matrix preserves the weight decomposition.
   
Let $\la_1,\la_2\not\in \La^+$. Let
\be
V^q_{\la_1}\T
V^q_{\la_2}=\bigoplus\limits_{l=0}^{\infty} V^q_{\la_1+\la_2-l}
\ee
be the decomposition of the tensor product of $U_qsl_2$ Verma modules into the
direct sum of irreducible modules, and let
$\Pi^{(l)}$ be the projector on $V^q_{\la_1+\la_2-l}$ along the other summands.

There is a formula for the trigonometric $R$-matrix:
\be
R^q_{V^q_{\la_1}V^q_{\la_2}}(x)=R^q_{\la_1,\la_2}(0)
\sum_{l=0}^\infty\Pi^{(l)}\prod_{s=0}^{l-1}
\frac{1-xq^{2s-2\la_1-2\la_2}}{1-xq^{2\la_1+2\la_2-2s}},
\ee
where
\be
R^q_{\la_1,\la_2}(0)=
q^{2\la_1\la_2-2h\T h}
\sum_{k=0}^\infty \, (q^2-1)^{2k}\,
\prod_{s=1}^k(1-q^{2s})^{-1}\,
(q^hf_q\T q^{-h}e_q)^k,
\ee
see \cite{CP},\cite{T}.

Let $\la_1,\la_2\in \La^+$. Let
\be
L^q_{\la_1}\T L^q_{\la_2}=\bigoplus\limits_{l=0}^{2\min\{\la_1,\la_2\}}L^q_{\la_1+\la_2-l}
\ee
be the decomposition of the tensor product of finite dimensional irreducible $U_qsl_2$ modules into
the direct sum of irreducibles, and let
$\Pi^{(l)}$ be the projector on $L^q_{\la_1+\la_2-l}$ along the other summands.

There is a formula for the trigonometric $R$-matrix,
\be
R^q_{L^q_{\la_1}L^q_{\la_2}}(x)=R^q_{\la_1,\la_2}(0)
\sum_{l=0}^{2\min\{\la_1,\la_2\}}\Pi^{(l)}\prod_{s=0}^{l-1}
\frac{1-xq^{2s-2\la_1-2\la_2}}{1-xq^{2\la_1+2\la_2-2s}},
\ee
see \cite{CP},\cite{T}.
The vector spaces $V^q_{\la_1}\T V^q_{\la_2}$ for different values of $\la_1,\la_2$ are identified by
distinguished bases $\{f^{l_1}v^q_1\T f^{l_2}v^q_2\,|\, l_1,l_2\in\Z_{\ge 0}\}$.
   
\begin{thm}\label{trigonometric $R$-matrix}
Let $q\in\C$ be not a root of unity.

1. The trigonometric $R$-matrix $R^q_{V^q_{\la_1}V^q_{\la_2}}(x)\in \End(V\T V)$ is
a meromorphic
function of $x,\la_1,\la_2$. The poles of
$R^q_{V^q_{\la_1}V^q_{\la_2}}(x)$  have the form
$x=q^{-2\la_1-2\la_2+2s}$, where $s\in\Z_{\ge 0}$.

2. Let $x$ be generic. Then
the trigonometric $R$-matrix $R^q_{V^q_{\la_1}V^q_{\la_2}}(x)$ preserves
$S^q_{\la_1}\T V^q_{\la_2}+V^q_{\la_1}\T S^q_{\la_2}$.

3. Let $V^q_{\la_1}\T V^q_{\la_2}\to L^q_{\la_1}\T L^q_{\la_2}$ be the canonical factorization map.
Then for generic $x$, the trigonometric $R$-matrix  $R^q_{V^q_{\la_1}V^q_{\la_2}}(x)$ can be
factorized to an operator
$R^q(x)\in \End(L^q_{\la_1}\T L^q_{\la_2})$ and, moreover,
$R^q(x)=R^q_{L^q_{\la_1}L^q_{\la_2}}(x)$.
\end{thm}
An analog of Theorem~\ref{rational $R$-matrix} for the elliptic $R$-matrix is proved in Theorems 8, 31
in \cite{FTV1}.
The same proof works for the trigonometric $R$-matrix. 
 
\subsection{The qKZ equation.}\label{qKZ}

The \emph {rational quantized Knizhnik-Zamolodchikov equation(qKZ)} associated 
to $sl_2$ is the following system of linear
difference equations for a
function $\Psi(z_1,\dots,z_n)$ with values in a tensor product 
$M_1\T\ldots\T M_n$ of $sl_2$ modules: 
\bea
\Psi(z_1,\dots,z_m+p,\dots,z_n)=R_{M_m,M_{m-1}}(z_m-z_{m-1}+p)\dots 
R_{M_m,M_1}(z_m-z_1+p)e^{-\mu h_m}\times 
\\
\times R_{M_m,M_n}(z_m-z_n)\dots R_{M_m,M_{m+1}}(z_m-z_{m+1})\Psi(z_1,\dots,z_n),
\eea
for $m=1,\dots,n$. Here $p$, $\mu$ are complex parameters, $\mu$ is chosen so that $0\le \Imag\mu 
<2\pi$;
$h_m$ is the operator $h\in sl_2$ acting in 
the $m$-th representation, 
$R_{M_iM_j}(x)\in \End(M_i\T M_j)$ is the rational $R$-matrix acting in the $i$-th and $j$-th factors, 
see \cite{FR}.
The linear operators in the right hand side of the equations are called the \emph{qKZ operators}.

The qKZ operators commute with the action of the operator $h\in sl_2$ in the tensor product
$M_1\T\ldots\T M_n$.
Therefore, in order to construct all solutions of the qKZ equation,
it is enough to solve the qKZ equation with values in weight spaces $(M_1\T\ldots\T M_n)_l$. 

If the parameter
$\mu$ of the equation is equal to zero, then the qKZ operators commute with 
the $sl_2$ action in the tensor product $M_1\T\ldots\T M_n$,
and in order to construct all solutions of the qKZ equation
in this case, it is enough to solve the equation with values in 
singular weight spaces $(M_1\T\ldots\T M_n)_l^{sing}$.

Let $\pi:\,V_{\la_1}\T\ldots\T V_{\la_n}\to L_{\la_1}\T\ldots\T L_{\la_n}$
be the canonical projection map.

\begin{lemma}\label{project}
Let $\Psi(z)$ be a solution of the qKZ equation with values in $V_{\la_1}\T\ldots\T V_{\la_n}$. Then
$\pi\circ\Psi(z)$ is a solution of the qKZ equation with values in
$L_{\la_1}\T\ldots\T L_{\la_n}$.
\end{lemma}

Lemma~\ref{project} follows from Theorem~\ref{rational $R$-matrix}.

\section{The Hypergeometric pairing, \cite{TV1} }\label{pairing}

\subsection{The phase function}

Let $z=(z_1,\dots,z_n)\in\C^n,\;\la=(\la_1,\dots,\la_n)\in\C^n ,\; 
t=(t_1,\dots,t_l)\in\C^l$.
The \emph{phase function} is defined by the following formula:
\be
\Phi_l(t,z,\la)=
\exp(\mu\sum_{i=1}^lt_i/p)\prod_{i=1}^n\prod_{j=1}^l
\frac{\Gamma((t_j-z_i+\la_i)/p)}{\Gamma((t_j-z_i-\la_i)/p)}
\prod_{1\le i< j\le l}\frac{\Gamma((t_i-t_j-1)/p)}{\Gamma((t_i-t_j+1)/p)}.
\ee

\subsection {Actions of the symmetric group}\label{actions}

Let $f=f(t_1,\dots,t_l)$ be a function. For a permutation $\sigma\in{\Bbb S}^l$, define the functions
$[f]_{\sigma}^{rat}$ and $[f]_{\sigma}^{trig}$ via the action of the simple transpositions 
$(i,i+1)\in {\Bbb S^l}$, $i=1,\dots,\l-1$, given by      

\be
[f]_{(i,i+1)}^{rat}(t_1,\dots, t_l)=
f(t_1,\dots,t_{i+1},t_i,\dots,t_l)\frac{t_i-t_{i+1}-1}{t_i-t_{i+1}+1},
\ee
\be
[f]_{(i,i+1)}^{trig}(t_1,\dots, t_l)=
f(t_1,\dots,t_{i+1},t_i,\dots,t_l)
\frac{
\sin(\pi(t_i-t_{i+1}-1)/p)}{\sin(\pi(t_i-t_{i+1}+1)/p)}.
\ee

If for all $\sigma\in{\Bbb S}^l$, $[f]_{\sigma}^{rat}=f$, we will say that the function is
\emph{symmetric with respect to the rational action}.
If for all $\sigma\in{\Bbb S}^l$, $[f]_{\sigma}^{trig}=f$, we will say that the function is
\emph{symmetric with respect to the trigonometric action}.
 
This definition implies the following important Remark.

{\bf Remark.} If $w(t_1,\dots, t_l)$ is ${\Bbb S}^l$ symmetric 
with respect to the 
rational action and $W(t_1,\dots, t_l)$ is ${\Bbb S}^l$ symmetric with respect to the 
trigonometric action,
then $\Phi_lwW$ is a symmetric function  of $t_1,\dots, t_l$ (in the 
usual sense).

\subsection {Rational weight functions}

Fix natural numbers $n,l$.

Set $\Zb^n_l=\{\bar{l}=(l_1,\dots,\l_n)\in\Z^n_{\ge 0}\,|\, 
\sum\limits_{i=1}^nl_i=l\}$. For $\bar{l}\in\Zb^n_l$ 
and $m=0,1,\dots,n$, set $l^m=\sum\limits_{i=1}^ml_i$.

For $\bar{l}\in\Zb_l^n$, define the \emph{rational weight function} 
$w_{\bar{l}}$ by 
\bea
w_{\bar{l}}(t,z,\la)=  
\sum_{\sigma \in {\Bbb S}^l}\left[\prod_{m=1}^n\frac{1}{l_m!}\prod_{j=l^{m-1}+1}^{l^m}
\left( \frac{1}{t_j-z_m-\la_m}\prod_{k=1}^m\frac{t_j-z_k+\la_k}{t_j-z_k-\la_k}
\right) \right] _\sigma^{rat}.
\eea

For fixed $z,\la\in\C^n$, the space spanned over $\C$
by  all rational weight functions 
$w_{\bar{l}}(t,z,\la)$, $\bar{l}\in\Zb_l^n$, is called 
the \emph{hypergeometric rational space specialized at $z,\la$} and 
is denoted $\F(z,\la)=\F^n_l(z,\la)$.
This space is a space of functions of variable $t$.

For generic values of $z,\la$, the space $\F(z,\la)$ can be 
identified with the space $(V_{\la_1}^*\T\ldots\T V_{\la_n}^*)_l$ by the map 
\be
{\frak a}(z,\la):\; (f^{l_1}v_1)^*\T\ldots\T (f^{l_n}v_n)^* \mapsto w_{\bar{l}}(t,z,\la).
\ee
Here $\{(f^lv_i)^*\,|\,l\in\Z_{\ge 0}\}$ is the basis of $V_{\la_i}^*$, dual to the standard 
basis of $V_{\la_i}$ given by $\{f^lv_i\,|\,l\in\Z_{\ge 0}\}$,
see Lemma 4.5, Corallary 4.8 in \cite{TV1}.

\subsection{Trigonometric weight functions}\label{trig}
Fix natural numbers $n, l$.

For $\bar{l}\in\Zb^n_l$, 
define the \emph{trigonometric weight function} $W_{\bar{l}}$ by
\bea
\lefteqn{W_{\bar{l}}(t,z,\la)=}
\\
&&
\sum_{\sigma \in {\Bbb S}^l}
\left[\prod_{m=1}^n\prod_{s=1}^{l_m}
\frac{\sin(\pi/p)}{\sin(\pi s/p)}
\prod_{j=l^{m-1}+1}^{l^m}\!\!
\frac{\exp(\pi i(z_m-t_j)/p)}{\sin(\pi(t_j-z_m-\la_m)/p)}
\prod_{k=1}^m\frac{\sin(\pi(t_j-z_k+\la_k)/p)}{\sin(\pi(t_j-z_k-\la_k)/p)}
\right] _\sigma^{trig}.
\eea

A function $W(t,z,\la)$ is said to be a \emph{holomorphic trigonometric weight function} if
\bean\label{trig.decomposition}
W(t,z,\la)=\sum_{m_1+\ldots +m_n=l}a_{\bar{m}}(\la,e^{2\pi iz_1/p},\dots,e^{2\pi iz_n/p})
W_{\bar{m}}(t,z,\la),
\eean
where $a_{\bar{m}}(\la,u)$ are holomorphic functions of parameters $\la,u\in\C^n$.
We denote $\G$ the space of all holomorphic trigonometric weight functions. This space
is a space of functions of variables of $t,z,\la$.

For a permutation $\sigma\in{\Bbb S}^n$,
define the \emph{trigonometric weight functions}
$W^\sigma_{\bar{l}}$ by $W^\sigma_{\bar{l}}(t,z,\la)=W_{\sigma\bar{l}}(t,\sigma z,\sigma\la)$, where
$\sigma\bar{l}=(l_{\sigma_1},\dots,l_{\sigma_n})$,
$\sigma z=(z_{\sigma_1},\dots,z_{\sigma_n})$,
$\sigma\la=(\la_{\sigma_1},\dots,\la_{\sigma_n})$. 

For fixed $\la, z\in\C^n$ and $\sigma\in{\Bbb S}^n$, the space spanned over $\C$
by  all trigonometric weight functions
$W^\sigma_{\bar{l}}(t,z,\la),\,\bar{l}\in\Zb_l^n$, is called
the \emph{hypergeometric trigonometric space specialized at $z,\la$} and is denoted
$\G^\sigma(z,\la)=\G^{n,\sigma}_l(z,\la)$. This space is a space of functions of variable $t$.

For generic $z,\la$, the space $\G^\sigma (z,\la)$ does not depend on $\sigma$, see Section 2 in
\cite{TV1}, we denote it $\G (z,\la)$.
 
For $\bar{l}\in\Zb^{n-1}_l$, define the \emph{singular trigonometric weight function}
$W_{\bar{l}}^{sing}$ by
\bea
\lefteqn{W_{\bar{l}}^{sing}(t,z,\la)=}
\\
&&
\sum_{\sigma \in {\Bbb S}^l}
\left[\prod_{m=1}^{n-1}\prod_{s=1}^{l_m}
\frac{\sin(\pi/p)}{\sin(\pi s/p}
\sin(\pi (z_m-\la_m-z_{m+1}-\la_{m+1}+s-1)/p)\times\right.
\\&&
\times
\prod_{j=l^{m-1}+1}^{l^m}\!\!
\frac{1}{\sin(\pi(t_j-z_m-\la_m)/p)\sin(\pi(t_j-z_{m+1}-\la_{m+1})/p)}\times
\\&&
\left.\times 
\prod_{k=1}^m\frac{\sin(\pi(t_j-z_k+\la_k)/p)}{\sin(\pi(t_j-z_k-\la_k)/p)}
\right] _\sigma^{trig}.
\eea

For fixed $z,\la\in\C^n$, the space spanned over $\C$ by all singular
trigonometric weight
functions $W_{\bar{l}}(t,z,\la),\,\bar{l}\in\Zb_l^n$,
is called the \emph{singular hypergeometric trigonometric space specialized at $z,\la$} and
is denoted $\G^{sing}(z,\la)=\G^{sing,n}_l(z,\la)$.

We have $\G^{sing}(z,\la)\subset\G(z,\la)$, see Lemma 2.29 in
\cite{TV1}.

A function $W(t,z,\la)$ is said to be a \emph{holomorphic singular trigonometric weight
function} if
$W(t,z,\la)$ is a holomorphic trigonometric weight function,
and for all $z,\la\in\C^n$, the function $W(t,z,\la)$ belongs to $\G^{sing}(z,\la)$.
We denote $\G^{sing}$ the space of all holomorphic singular trigonometric weight functions. This
space is a space of functions of variables $t,z,\la$.

\begin{lemma}\label{W^{sing}}
For any
$\bar{l}\in\Zb^{n-1}_l$, the function
$W_{\bar{l}}^{sing}(t,z,\la)$ belongs to $\G^{sing}$.
\end{lemma}
\begin{proof}
By Lemmas 2.28 and 2.29 of \cite{TV1}, we have the decomposition \Ref{trig.decomposition},
\be
W_{\bar{l}}^{sing}(t,z,\la)=\sum_{m_1+\ldots +m_n=l}a_{\bar{m}}(\la,e^{2\pi iz_1/p},\dots,e^{2\pi iz_n/p})
W_{\bar{m}}(t,z,\la).
\ee  
Here $a_{\bar{m}}(\la,u)$ are meromorphic functions of $\la,u\in\C^n$. The functions $a_{\bar{m}}(\la,u)$
are holomorphic for the following reason.

Suppose, for some $\bar{m}_0\in\Zb^n_l$ and
some
$Z,\La\in\C^n$ the function  $a_{\bar{m}}(\la,u)$ has a pole at $Z,\La$. The functions
$W_{\bar{m}}(t,z,\la)$, $\bar{m}\in\Zb^n_l$ are
linearly independent for generic $t,z,\la$ by  Lemma 2.28 in \cite{TV1}. Moreover, for generic $t$, the
functions
$W_{\bar{m}}(t,z,\la)$, $\bar{m}\in\Zb^n_l$, do not have poles at $Z,\La$. Hence, for generic $t$, the
function
$W_{\bar{l}}^{sing}(t,z,\la)$ has a pole at $Z,\La$. But this is not so.
\end{proof}
For  $\bar{l}\in\Zb^{n}_l$, define the \emph{weight coefficient} $c_{\bar{l}}(\la)$ by
\be
c_{\bar{l}}(\la)=
\prod_{m=1}^n\prod_{s=0}^{l_m-1}
\frac{\sin(\pi(s+1)/p)\sin(\pi(2\la_m-s)/p)}{\sin(\pi/p)}.
\ee

Let $q=e^{\pi i/p}$.
The weight  space $(V_{\la_{\sigma_1}}^q\T\ldots\T V_{\la_{\sigma_n}}^q)_l$
is mapped to the 
trigonometric hypergeometric space $\G^\sigma(z,\la)$ 
by 
\bean\label{trig.ident.}
{\frak b}_\sigma(z,\la) :\; f_q^{l_{\sigma_1}}v_{\sigma_1}^q\T\ldots\T
f_q^{l_{\sigma_n}}v_{\sigma_n}^q \mapsto 
c_{\bar{l}}(\la)W^\sigma_{\bar{l}}(t,z,\la).
\eean
If  $z,\la$ are generic and $q$ is not a root of unity, then the
map ${\frak b}_\sigma(z,\la)$ is an isomorphism of vector spaces, see Lemma 4.17 and Corallary
4.20 in \cite{TV1}. Moreover, the singular hypergeometric space is
identified with the subspace of
singular
vectors,
\be
\G^{sing}(z,\la)={\frak b}_\sigma(z,\la)((V_{\la_{\sigma_1}}^q\T\ldots\T
V_{\la_{\sigma_n}}^q)_l^{sing}),
\ee
see Corallary 4.21 in \cite{TV1}.

The composition maps 
\be
{\frak b}_{\sigma,\sigma^\prime}(z,\la):\,
(V_{\la_{\sigma^\prime_1}}^q\T\ldots\T V_{\la_{\sigma^\prime_n}}^q)_l\to
(V_{\la_{\sigma_1}}^q\T\ldots\T V_{\la_{\sigma_n}}^q)_l,
\;\;
{\frak b}_{\sigma,\sigma^\prime}(z,\la)=({\frak b}_\sigma(z,\la))^{-1}
\circ{\frak b}_{\sigma^\prime}(z,\la),
\ee
are called the \emph{transition functions}.
\begin{thm}\label{cite1}
(Theorem 4.22 in \cite{TV1}.) Let $q$ be not a root of unity. 
For any $\sigma\in{\Bbb S}^n$ and any transposition $(m,m+1), m=1,\dots,n-1$, the transition function
\be
{\frak b}_{\sigma,\sigma\circ (m,m+1)}(z,\la):\,
(V_{\la_{\sigma_1}}^q\T\ldots\T V^q_{\la_{\sigma_{m+1}}}\T V_{\la_{\sigma_m}}^q\T\ldots\T
V_{\la_{\sigma_n}}^q)_l
\to (V_{\la_{\sigma_1}}^q\T\ldots\T V_{\la_{\sigma_n}}^q)_l
\ee 
equals the operator
$P_{V_{\la_{\sigma_{m+1}}}^qV_{\la_{\sigma_m}}^q}R_{V_{\la_{\sigma_{m+1}}}^qV_{\la_{\sigma_m}}^q}
(e^{2\pi i(z_{\sigma_{m+1}}-z_{\sigma_m})/p})$ acting in the $m$-th and $(m+1)$-st factors.
Here $P_{V_{\la_{\sigma_{m+1}}}^qV_{\la_{\sigma_m}}^q}$ is the operator of permutation of the
factors $V_{\la_{\sigma_{m+1}}}^q$ and $V_{\la_{\sigma_m}}^q$.
\end{thm}

\begin{corollary}\label{transition functions}
Let $q$ be not a root of unity.
For any $\sigma\in{\Bbb S}^n$ and any transposition $(m,m+1), m=1,\dots,n-1$,
the transition function ${\frak b}_{\sigma,\sigma^\prime}(z,\la)$ can be factorized to a
transition function
\be
{\frak B}_{\sigma,\sigma^\prime}(z,\la):\,
(L_{\la_{\sigma^\prime_1}}^q\T\ldots\T L_{\la_{\sigma^\prime_n}}^q)_l\to
(L_{\la_{\sigma_1}}^q\T\ldots\T L_{\la_{\sigma_n}}^q)_l.
\ee
Moreover,
for any $\sigma\in{\Bbb S}^n$ and any transposition $(m,m+1), m=1,\dots,n-1$, the induced transition
function
\be
{\frak B}_{\sigma,\sigma\circ(m,m+1)}(z,\la):\,
(L_{\la_{\sigma_1}}^q\T\ldots\T L^q_{\la_{\sigma_{m+1}}}\T L_{\la_{\sigma_m}}^q\T\ldots\T
L_{\la_{\sigma_n}}^q)_l
\to (L_{\la_{\sigma_1}}^q\T\ldots\T L_{\la_{\sigma_n}}^q)_l
\ee
equals the operator
$P_{L_{\la_{\sigma_{m+1}}}^qL_{\la_{\sigma_m}}^q}R_{L_{\la_{\sigma_{m+1}}}^qL_{\la_{\sigma_m}}^q}
(e^{2\pi i(z_{\sigma_{m+1}}-z_{\sigma_m})/p})$ acting in the $m$-th and $(m+1)$-th factors.
\end{corollary}

Corollary~\ref{transition functions} follows from Theorems~\ref{trigonometric $R$-matrix}
and \ref{cite1}.

Thus, the map
${\frak b}_\sigma(z,\la): \,(V_{\la_{\sigma_1}}^q\T\ldots\T V_{\la_{\sigma_n}}^q)_l\to\G(z,\la)$ 
can be factorized to a map
\bean\label{trig.ident.1}
{\frak B}_\sigma(z,\la): \,(L_{\la_{\sigma_1}}^q\T\ldots\T
L_{\la_{\sigma_n}}^q)_l\to\G(z,\la),       
\eean
defined by the same formula \Ref{trig.ident.}. Moreover, the image of ${\frak B}_\sigma(z,\la)$ does not
depend on $\sigma$. Transition functions 
$({\frak B}_\sigma(z,\la))^{-1}\circ{\frak B}_{\sigma^\prime}(z,\la)$ are well defined, and
$({\frak B}_\sigma(z,\la))^{-1}\circ{\frak B}_{\sigma^\prime}(z,\la)=
{\frak B}_{\sigma,\sigma^\prime}(z,\la)$.

\subsection{Hypergeometric integrals}\label{integrals}

Fix $p\in\C,\,\Real p<0$. Let $\Imag\mu\neq 0$.
Assume that the parameters $z,\la\in\C^n$ satisfy the condition $\Real (z_i+\la_i)<0$ and
$\Real (z_i-\la_i)>0$ for all $i=1,\dots,n$. For a 
rational weight function $w_{\bar{l}}(t,z,\la)$, $\bar{l}\in\Zb^n_l$, and a trigonometric
weight
function 
$W(t,z,\la)\in\G$, define the \emph{hypergeometric integral} $I(w,W)(z,\la)$ by the formula
\bean\label{int}
I(w,W)(z,\la)=\int\limits_{\Real t_i=0,\atop
i=1,\dots,l}\Phi_l(t,z,\la)w(t,z,\la)W(t,z,\la)\,d^lt,
\eean
where $d^lt=dt_1\ldots dt_l$.

The hypergeometric integral for generic $z,\la$ and an arbitrary step $p$ with negative real part is
defined by analytic continuation with respect to $z,\la$ and $p$. This analytic continuation makes
sense since the integrand is meromorphic in $z,\la$ and $p$. The poles of the integrand are located at 
the union of hyperplanes
\bean\label{hyperplanes}
t_i=z_k\pm(\la_k+sp),\qquad t_i=t_j\pm(1-sp),
\eean
$i,j=1,\dots,l,\,k=1,\dots,n,\,s\in\Z_{\ge 0}$. We move the parameters $z,\la$ and $p$ in such a way 
that the topology of the complement in $\C^l$ 
to the union of hyperplanes \Ref{hyperplanes} does not change. We deform
accordingly the integration cycle (the imaginary subspace) in such a way that it does not
intersect the hyperplanes
\Ref{hyperplanes} at any moment of the deformation. 
The deformed integration cycle is called the \emph{deformed imaginary 
subspace} and is denoted $\I (z,\la)$, $\, \I (z,\la) \subset \C^l$.
Then  the analytic continuation of integral \Ref{int} is given by
\bean\label{int'}
I(w,W)(z,\la)=\int\limits_{\I (z,\la)}\Phi_l(t,z,\la)w(t,z,\la)W(t,z,\la)\,d^lt,
\eean
see Section 5 in \cite{TV1}.

\begin{thm}\label{cite2}
(Theorem 5.7 in \cite{TV1}.)
Let $\Imag\mu\neq 0$. For any $w_{\bar{l}}(t,z,\la)$, $\bar{l}\in\Zb^n_l$, and any
$W(t,z,\la)\in\G$,
the hypergeometric integral \Ref{int'} is a univalued meromorphic function 
of variables $p, z, \la$ holomorphic 
on the set
\be
\Real p<0,\qquad \{1,\dots,l\}\not\subset p\Z,
\ee
\bean\label{old}
2\la_m-s\not\in p\Z,\qquad m=1,\dots,n,\qquad s=1-l,\dots,l-1,
\eean
\be
z_k\pm\la_k-z_m\pm\la_m-s\not \in p\Z,\qquad k,m=1,\dots,n, \qquad k\neq m,
\qquad s=1-l,\dots,l-1,
\ee
where we allow arbitrary combinations of $\pm$.
\end{thm}

The case $\Imag\mu=0$ is treated in the same manner. 
\begin{thm}\label{cite2'}
(Theorem 5.8 in \cite{TV1}.)
Let $\Imag\mu=0$. For any $w_{\bar{l}}(t,z,\la)$, $\bar{l}\in\Zb^n_l$, and any
$W(t,z,\la)\in\G^{sing}$,
the hypergeometric integral \Ref{int'} is a univalued meromorphic function of variables
$p, z, \la$ holomorphic
on the set \Ref{old}.
\end{thm}

For a function $W(t,z,\la)\in\G$, 
let $\Psi_W(z,\la)$ be the following $V_{\la_1}\T\ldots\T V_{\la_n}$-valued 
function
\bean\label{Psi}
\Psi_W(z,\la) \, = \, \sum_{l_1+\ldots+l_n=l} \,
I(w_{\bar{l}},W)(z,\la) \, f^{l_1}v_1\T\ldots\T f^{l_n}v_n\, .
\eean

\begin{thm}\label{cite3}
(Corollaries 5.25, 5.26 in \cite{TV1}.)
Let $p,z,\la$ be such that conditions \Ref{old} are satisfied.

i) Let $\Imag\mu\neq 0$. Then for any function 
$W\in\G$, the function $\Psi_W(z,\la)$ is a solution of the qKZ equation with values in 
$(V_{\la_1}\T\ldots\T V_{\la_n})_l$.

ii) Let $\Imag\mu=0$. Then for any function $W\in\G^{sing}$, the function
$\Psi_W(z,\la)$ is a solution
of the qKZ equation with values in 
$(V_{\la_1}\T\ldots\T V_{\la_n})_l^{sing}$. 
\end{thm}

The solutions of the qKZ equation defined by \Ref{Psi} are called the \emph{hypergeometric
solutions}.

\section{Main results}\label{main}

Fix natural numbers $n,l$.
We will often assume the following restrictions on the parameters $p,z,\la$:
\bean\label{step}
\Real p<0,\qquad 1\not\in p\Z,
\eean
\bean\label{weights1}
\{s\, | \, s\in\Z_{>0}, s<
2\max\{\Real\la_1,\dots,\Real\la_n\},\, s\le l\}\bigcap\{p\Z\}=\emptyset,
\eean
\bean\label{weights2}
\{2\la_m-s\,|\,s\in\Z_{\ge 0},\; s<2\Real\la_m,\, s<l\}\bigcap \{p\Z\}=\emptyset, \qquad
m=1,\dots,n,
\eean
\bean\label{resonance}
z_k-z_m\pm(\la_k+\la_m)+s \not\in\{p\Z\},\qquad k,m=1,\dots,n,\,k\neq m,\, s=1-l,\dots,l-1.
\eean

Sometimes, in addition we will make the following assumption. Let $\la\in\C^n$. For
each
$i\in\{1,\dots,n\}$ such that $\la_i\not\in\La^+$, assume
\bean\label{step3}
\{1,\dots,l\}\bigcap\{p\Z\}=\emptyset,
\eean
\bean\label{weights3}
\{2\la_i-s\,|\,s=0,1,\dots,l-1\}\bigcap \{p\Z\}=\emptyset.
\eean

Notice that if the parameters $p, z, \la$ satisfy conditions
\Ref{step}-\Ref{weights3},
then they also satisfy conditions \Ref{old}.

\subsection{Analytic continuation. The case $\Imag\mu\neq 0$}\label{an}

Let $\La=(\La_1,\dots,\La_n)\in\C^n$, $\bar{l}=(l_1,\dots,l_n)\in\Z^n_{\ge 0}$.
An $i$-th coordinate of $\bar{l}$ is called
\emph{$\La$-admissible} if either $\La_i\not\in\La^+$
or $\La_i\in\La^+$ and $l_i\le 2\La_i$.
The index $\bar{l}$ is called \emph{$\La$-admissible} if all its coordinates are
$\La$-admissible. Denote $B_\La(\bar{l})\subset\{1,\dots,n\}$ the set of
all non-$\La$-admissible coordinates of $\bar{l}$.

{\bf Remark.}
Let $p$ be generic, $\La\in\C^n,\,\bar{l}\in\Zb^n_l$. Then the weight 
coefficient $c_{\bar{l}}(\La)$ is not equal to zero if and only if the index
$\bar{l}$ is $\La$-admissible.

Let $\bar{m}\in\Zb^n_l$ and $B_\La(\bar{l})\subseteq B_\La(\bar{m})$. Then the function
$c(\la)=c_{\bar{m}}(\la)/c_{\bar{l}}(\la)$ is holomorphic at $\La$. Moreover,
$c(\La)=0$ if $B_\La(\bar{l})\neq B_\La(\bar{m})$, and $c(\La)$ is not zero if 
$B_\La(\bar{l})=B_\La(\bar{m})$.

\begin{thm}\label{an.cont.1}
Let $\Imag\mu\neq 0$. 
Let $p$ satisfy condition
 \Ref{step}. 
Let $Z, \La \in \C^n$ satisfy conditions \Ref{weights1}-\Ref{resonance}. 
Let  $\bar{l}, \bar{m}\in \Zb^n_l$.
Assume that for all $i=1,\dots,n$ either the $i$-th coordinate of $\bar{l}$ or 
the $i$-th coordinate of $\bar{m}$ is $\La$-admissible, i.e.
$B(\bar{l})\bigcap B(\bar{m})=\emptyset$.
Then the hypergeometric integral $I(w_{\bar{l}},W_{\bar{m}})(z,\la)$
is holomorphic at $Z, \La$.
Moreover, there exists a contour of integration $\I(Z,\La)\subset\C^l$ independent on
$\bar{l},\bar{m}$, such
that for all 
$z, \la$, in a small neighborhood of $Z, \La$ we
have
\be
I(w_{\bar{l}},
W_{\bar{m}})(z,\la)=\int\limits_{\I(Z,\La)}\Phi_l(t,z,\la)w_{\bar{l}}(t,z,\la)W_{\bar{m}}(t,z,\la)\,d^lt.
\ee
\end{thm}

A contour $\I(Z,\La)$ with properties indicated in Theorem~\ref{an.cont.1} is called an
\emph{integration contour associated to $Z,\La$}.

Theorem~\ref{an.cont.1} is proved in Section~\ref{proofs}.

For $\bar{l},\bar{m}\in \Zb^n_l$, introduce a function
\bean\label{c int}
J_{\bar{l},\bar{m}}(z,\la)=c_{\bar{m}}(\la)I(w_{\bar{l}},W_{\bar{m}})(z,\la).
\eean

\begin{thm}\label{an.cont.2}
Let $\Imag\mu\neq 0$. 
Let $p$ satisfy condition
 \Ref{step}.
Let $Z, \La \in \C^n$ satisfy conditions \Ref{weights1}-\Ref{resonance}.
Let  $\bar{l}, \bar{m}\in \Zb^n_l$.
Then the meromorphic function $J_{\bar{l},\bar{m}}(z,\la)$
is holomorphic at $Z, \La$.
Moreover, if $B_\La(\bar{l})\subset B_\La(\bar{m})$ and $B_\La(\bar{l})\neq B_\La(\bar{m})$, then 
$J_{\bar{l},\bar{m}}(Z,\La)=0$.
\end{thm}

Theorem~\ref{an.cont.2} is proved in Section~\ref{proofs}.

Let $\La\in\C^n,\,\bar{l},\bar{m}\in \Zb^n_l$ and $B_\La(\bar{l})=B_\La(\bar{m})$. 
Denote $B$
the
set $B_\La(\bar{l})=B_\La(\bar{m})$.
Introduce $\La^\prime(B)=(\La_1^\prime,\dots,\La_n^\prime)\in\C^n$ by the rule: 
$\La_i^\prime=\La_i$ if $i\not\in B$ and $\La_i^\prime=-\La_i-1$ if
$i\in B$.

Let $\bar{l}^\prime(B)=(l_1^\prime,\dots,l_n^\prime)$, where
$\l_i^\prime=l_i$ if $i\not\in B$ and $l_i^\prime=l_i-2\La_i-1$ if
$i\in B$. Similarly, let $\bar{m}^\prime(B)=(m_1^\prime,\dots,m_n^\prime)$, where
$m_i^\prime=m_i$ if $i\not\in B$ and $m_i^\prime=m_i-2\La_i-1$ if 
$i\in B$.

We have 
$\l_1^\prime+\ldots +l_n^\prime=m_1^\prime+\ldots +m_n^\prime=
l-2\sum\limits_{i\in B}\La_i-|B|$. We denote this number $l^\prime(B)$.

\begin{thm}\label{functor}
Let $\Imag\mu\neq 0$.
Let $p$ satisfy condition
 \Ref{step}.
Let $Z, \La \in \C^n$ satisfy conditions \Ref{weights1}-\Ref{resonance}.
Let  $\bar{l}, \bar{m}\in \Zb^n_l$ and $B_\La(\bar{l})=B_\La(\bar{m})=B$. Then
\be
J_{\bar{l},\bar{m}}(Z,\La)=
{\frak C}_\La(Z)
J_{\bar{l}^\prime(B),\bar{m}^\prime(B)}(Z,\La^\prime(B)),
\ee   
where ${\frak C}_\la (z)$ is a nonzero holomorphic function at $Z$ given below.
\end{thm}

Theorem~\ref{functor} is proved in Section~\ref{proofs}.

Notice that
$J_{\bar{l}^\prime(B),\bar{m}^\prime(B)}(z,\La^\prime(B))=c_{\bar{m}^\prime(B)}(\La^\prime(B))
I(w_{\bar{l}^\prime(B)},W_{\bar{m}^\prime(B)})(z,\La^\prime(B))$, where
\newline
$I(w_{\bar{l}^\prime(B)},W_{\bar{m}^\prime(B)})(z,\la)$ is
an $l^\prime(B)$-dimensional hypergeometric integral. The indices 
$\bar{l}^\prime(B)$, $\bar{m}^\prime(B)$  
are $\La^\prime(B)$-admissible. By Theorem~\ref{an.cont.1} the integral
$I(w_{\bar{l}^\prime(B)},W_{\bar{m}^\prime(B)})(z,\La^\prime(B))$ is well defined.
Theorem~\ref{functor} connects $l$- and $l^\prime(B)$-dimensional hypergeometric integrals.

Now, we describe the function ${\frak C}_\la(z)$. For $k\in\Z_{\ge 0}$, let
\be
\psi_k=\frac{-1}{\pi^{k+2}\,k!\,(k+1)!}\,
\Gamma(-(k+1)/p)\,\Gamma(1/p)^{k+1}\,
\prod_{j=1}^{k+1}\sin(j\pi/p).
\ee
For $z,\la\in\C^n,j\in\{1,\dots,n\},k\in\Z_{\ge 0}$, let
\be
\phi_{\la,j,k}(z)=\prod_{s=0}^{k}
\left(\prod_{i=0}^{j-1}
\frac{\Gamma((z_i-z_j+\la_i+\la_j-s)/p)}{\Gamma((z_i-z_j-\la_i-\la_j+s)/p)}
\prod_{i=j+1}^n
\frac{\Gamma((z_j-z_i+\la_j+\la_i-s)/p)}{\Gamma((z_j-z_i-\la_j-\la_i+s)/p)}\right).
\ee

Then 
\be
{\frak C}_\la(z)=
\frac{l!}{(l^\prime(B)) !}
\prod_{j\in B} 
e^{\mu(2\la_j+1)z_j\pi i /p}
\psi_{2\la_j}\,\phi_{\la,j,2\la_j}(z).
\ee

Consider a square matrix
\bean\label{matrix}
J^l(z,\la)=\{J_{\bar{l},\bar{m}}(z,\la)\}_{\bar{l},\bar{m}\in\Zb_l^n}.
\eean

Recall that for generic $p$, we have $c_{\bar{m}}(\la)\neq 0$ if $\bar{m}$ is $\la$-admissible.
According to
Theorem~\ref{cite2}, the matrix \Ref{matrix} is holomorphic if conditions \Ref{old} hold. According to
Theorem~\ref{an.cont.2}, the matrix \Ref{matrix} is holomorphic on a larger set of parameters
\Ref{step}-\Ref{resonance}.

According to Theorem~\ref{an.cont.2}, if some of
$\la_1,\dots\,\la_n$ become nonnegative half-integers, then the matrix
\Ref{matrix} becomes upper block triangular in the following sense. 

Let $\La\in\C^n$. Divide the set of indices $\Zb_l^n$ into subsets labeled by subsets of
$\{1,\dots,n\}$.
The subset of  $\Zb_l^n$ corresponding to a subset $B\subset\{1,\dots,n\}$ consists of all
$\bar{l}\in\Zb_l^n$ such that $B_\La(\bar{l})=B$.
Since the rows and columns of matrix \Ref{matrix} are labeled by elements of $\Zb_l^n$, matrix
\Ref{matrix} is divided into blocks labeled by pairs 
$B_1,B_2\subset\{1,\dots,n\}$. 

Choose any order $<$ on the set of all subsets of $\{1,\dots,n\}$ such that
for any $B_1,B_2\subset\{1,\dots,n\},B_1\subset B_2$, we have $B_1<B_2$.
Theorem~\ref{an.cont.2} says that matrix \Ref{matrix} is upper block triangular at
$(z,\La)$ with
respect to this order. Namely, the block corresponding to a pair $B_1,B_2$ is equal to zero if
$B_1<B_2$.

Theorem~\ref{functor} describes the diagonal blocks of matrix \Ref{matrix}.
The diagonal block of $J^l(z,\La)$ corresponding to a subset $B\subset\{1,\dots,n\}$ coincides 
(up to multiplication from the left by the non-degenerate diagonal matrix
$\{{\frak C}_{\La,\bar{m}}(z)\}_{\bar{m}\in\Zb_l^n}$) with the 
diagonal block of the matrix $J^{l^\prime(B)}(z,\La^\prime(B))$ corresponding to the empty subset of
$\{1,\dots,n\}$. 

\subsection{Analytic continuation}\label{an.'}
{\bf The case $\Imag\mu=0$} is treated similarly.

A function $W\in\G^{sing}$ is called $\La$-\emph{admissible} if for all non-$\La$-admissible 
$\bar{m}\in\Zb^n_l$, the functions
$a_{\bar{m}}(\la,u)$
in decomposition \Ref{trig.decomposition} are equal to zero.

\begin{thm}\label{an.cont.'1}
Let $\Imag\mu= 0$.
Let $p$ satisfy condition \Ref{step}.
Let $Z, \La \in \C^n$ satisfy conditions \Ref{weights1}-\Ref{resonance}.
Let  $\bar{l}\in\Zb^n_l$. Let $W\in\G^{sing}$ be $\La$-admissible.
Then the hypergeometric integral $I(w_{\bar{l}},W)(z,\la)$
is holomorphic at $Z, \La$.
Moreover, there exists a contour of integration $\I(Z,\La)\subset\C^l$ independent on $\bar{l}$ and
$W$, such
that for all
$z, \la$, in a small neighborhood of $Z,\La$ we
have
\be
I(w_{\bar{l}},W)(z,\la)=
\int\limits_{\I(Z,\La)}\Phi_l(t,z,\la)w_{\bar{l}}(t,z,\la)W(t,z,\la)\,d^lt.
\ee
\end{thm}

A contour $\I(Z,\La)$ with properties indicated in Theorem~\ref{an.cont.'1} is called an
\emph{integration contour associated to $Z,\La$}.

The proof of Theorem~\ref{an.cont.'1} is similar to the proof of Theorem~\ref{an.cont.1}.  

It follows from the proof that a contour of integration associated to $Z,\La$ with respect to
Theorem~\ref{an.cont.1} is also an integration contour associated to $Z,\La$ with respect to
Theorem~\ref{an.cont.'1} and vice versa.

\begin{thm}\label{an.cont.'2}
Let $\Imag\mu= 0$.
Let $p$ satisfy condition
 \Ref{step}.
Let $Z, \La \in \C^n$ satisfy conditions \Ref{weights1}-\Ref{resonance}.
Let  $\bar{l}\in\Zb^n_l$, $W\in\G^{sing}$.
Then the meromorphic function $c_{\bar{l}}(\la)I(w_{\bar{l}},W)(z,\la)$
is holomorphic at $Z, \La$.
Moreover, if $\bar{l}$ is not $\La$-admissible and 
$W$ is $\La$-admissible,
then the function
$c_{\bar{l}}(\la)I(w_{\bar{l}},W)(z,\la)$ is equal to zero at $(Z,\La)$.
\end{thm}

The proof of Theorem~\ref{an.cont.'2} is similar to the proof of Theorem~\ref{an.cont.2}.

Notice that in Theorem~\ref{an.cont.'2} we consider the function $c_{\bar{l}}(\la)I(w_{\bar{l}},W)(z,\la)$
and not the function $c_{\bar{m}}(\la)I(w_{\bar{l}},W_{\bar{m}})(z,\la)$ as in Theorem~\ref{an.cont.2}.

\subsection{Solutions of qKZ with values in irreducible 
representations.}\label{fin}

For $\bar{m}\in\Zb^n_l$, consider a $V_{\la_1}\T\ldots\T V_{\la_n}$-valued
function 
\bean\label{cPsi}
\Psi_{\bar{m}}(z,\la)=\sum_{l_1+\ldots+l_n=l} 
J_{\bar{l},\bar{m}}(z,\la) f^{l_1}v_1\T\ldots\T f^{l_n}v_n.
\eean

\begin{corollary}\label{cor}
Let $\Imag\mu\ne 0$.
Let $p\in\C$ and $\La\in\C^n$ satisfy conditions \Ref{step}-\Ref{weights2}.
Then the function
$\Psi_{\bar{m}}(z,\La)$ given by \Ref{cPsi}
is a meromorphic solution of the qKZ equation with values in $(V_{\La_1}\T\ldots\T V_{\La_n})_l$, holomorphic
for all $z$
satisfying condition \Ref{resonance}.
\end{corollary}

Corollary~\ref{cor} follows from
Theorems~\ref{rational $R$-matrix}, \ref{cite3}, and
\ref{an.cont.2}.

Let $\La\in\C^n$. For any $\La$-admissible $\bar{m}\in\Zb^n_l$, consider a function
\bean\label{Psi fin}
\Psi_{\bar{m}}(z,\La)=\sum
I(w_{\bar{l}},W_{\bar{m}})(z,\la) f^{l_1}v_1\T\ldots\T f^{l_n}v_n,
\eean
where the sum is over all $\La$-admissible $\bar{l}\in\Zb^n_l$. 

\begin{corollary}\label{cor1}
Let $\Imag\mu\ne 0$.
Let $p\in\C$ and $\La\in\C^n$ satisfy conditions \Ref{step}-\Ref{weights2}.
Then for any $\La$-admissible $\bar{m}\in\Zb^n_l$,
the function $\Psi_{\bar{m}}(z,\La)$ given by \Ref{Psi fin}
is a meromorphic solution of the qKZ equation with values in $(L_{\La_1}\T\ldots\T L_{\La_n})_l$, holomorphic
for all $z$
satisfying condition \Ref{resonance}.
\end{corollary}

Notice that vectors $f^{l_1}v_1\T\ldots\T f^{l_n}v_n$ with $\La$-admissible $\bar{l}\in\Zb^n_l$
form a basis in $(L_{\La_1}\T\ldots\T L_{\La_n})_l$.

Corollary~\ref{cor1} follows from Corollary~\ref{cor} and Lemma~\ref{project}.

\bigskip 

Let $\Imag\mu=0$. For $W\in\G^{sing}$, consider a $V_{\la_1}\T\ldots\T V_{\la_n}$-valued
function
$\Psi_W(z,\la)$ given by \Ref{Psi}.

\begin{corollary}\label{cor'}
Let $\Imag\mu=0$.
Let $p\in\C$ and $\La\in\C^n$ satisfy conditions \Ref{step}-\Ref{weights2}.   
Let $W\in\G^{sing}$ be $\La$-admissible.
Then the function $\Psi_W(z,\La)$ given by \Ref{Psi}
is a meromorpfic solution of the qKZ equation with values in $(V_{\La_1}\T\ldots\T V_{\La_n})^{sing}_l$,
holomorphic for
all $z$ satisfying condition \Ref{resonance}.
\end{corollary}

Corollary~\ref{cor'} follows from
Theorems~\ref{rational $R$-matrix},
\ref{cite3},
and
\ref{an.cont.'1}.

Let $\La\in\C^n,\,W\in\G^{sing}$. Consider a function
\bean\label{Psi fin'}
\Psi_W(z,\La)=\sum
I(w_{\bar{l}},W)(z,\la) f^{l_1}v_1\T\ldots\T f^{l_n}v_n,
\eean
where the sum is over $\La$-admissible $\bar{l}\in\Zb^n_l$.

\begin{corollary}\label{cor1'}
Let $\Imag\mu=0$.
Let $p\in\C$ and $\La\in\C^n$ satisfy conditions \Ref{step}-\Ref{weights2}.
Let $W\in\G^{sing}$ be $\La$-admissible.
Then $\Psi_W(z,\La)$ given by \Ref{Psi fin'} is a meromorphic
solution of the qKZ equation with values in $(L_{\La_1}\T\ldots\T L_{\La_n})^{sing}_l$, holomorphic
for all $z$
satisfying condition \Ref{resonance}.
\end{corollary}

Corollary~\ref{cor1'} follows from Corollary~\ref{cor'} and Lemma~\ref{project}.

The solutions of the qKZ equation defined by formulas \Ref{Psi fin},\Ref{Psi fin'} are called
the \emph{hypergeometric solutions}.

\subsection{Determinant of the hypergeometric pairing.}\label{determinant}

Let $\Imag\mu\ne 0$.
The determinant of the matrix $J^l(z,\la)$ defined by \Ref{matrix} 
is given by
\bea
\lefteqn
{\det(J_l(z,\la))=(2i)^{l{n+l-1 \choose n-1}}(l!)^{n+l-1 \choose n-1}
(e^\mu-1)^{-2\sum\limits_{m=1}^n\la_m/p\cdot{n+l-1\choose n}+
2n/p\cdot{n+l-1\choose n+1}}\times}
\\&&
\times \exp(\mu\sum_{m=1}^nz_m/p\cdot{n+l-1\choose n}\times
\\&&
\times \exp\left( (\mu+\pi
i)\left(\sum_{m=1}^n\la_m/p\cdot{n+l-1\choose n}-n/p\cdot{n+l-1\choose
n+1}\right)\right)\times
\\&&  
\times
\prod_{s=0}^{l-1}\left[ \Gamma(1+1/p)^n\Gamma(1+(s+1)/p)^{-n}   
\prod_{m=1}^n\frac{\pi}{\Gamma(1-(2\la_m-s)/p)}\right.\times
\\&&
\times\left.\prod_{1\le k<m\le n}
\frac{\Gamma((z_k+\la_k-z_m+\la_m-s)/p)}{\Gamma((z_k-\la_k-z_m-\la_m+s)/p)}
\right]^{n+l-s-2 \choose n-1},
\eea
cf. Theorem 5.14 of \cite{TV1}.

For $\La\in\C^n$, consider a matrix
\bean\label{adm.matrix}
J^l_{\rm adm}(z,\la)=\{J_{\bar{l},\bar{m}}(z,\la)\},
\eean
where $\bar{l},\bar{m}\in\Zb^n_l$ run through the set of all $\La$-admissible indices.

By Theorem~\ref{an.cont.1}, the matrix $J^l_{\rm adm}(z,\la)=J^l_{\rm adm}(z,\la,p)$ is holomorphic at
$(z,\La)$ if parameters $p,z,\La$ satisfy conditions \Ref{step}-\Ref{resonance}.

\begin{thm}\label{det}
Let $\Imag\mu\ne 0$. Let $p\in\C$ and $z,\La\in\C^n$ satisfy conditions
\Ref{step}-\Ref{weights3}.
Then the matrix $J^l_{\rm adm}(z,\La)$ is non-degenerate. Moreover, for generic $p$,
\bean\label{det-prod}
\det (J^l_{\rm adm}(z,\La))=
\prod\limits_{A,\,A\subset
B(\La)}\left(C_\La(A)(z)\det
(J^{l^\prime(A)}(z,\La^\prime(A)))\right)^{(-1)^{|A|}},
\eean
where $B(\La)$ is the set of all
$i\in\{1,\dots,n\}$ such that $\La_i\in\La^+$, and the function $C_\La(A)(z)$ is given by
\be
C_\La(A)(z)=
\prod_{\bar{m},\,B(\bar{m})\supseteq A}\,
\prod_{i\in A}\left(\frac{l!}{(l-2\La_i-1)!}\;\psi_{2\La_i}\,\phi_{\La,i,2\La_i}(z)\right).
\ee
\end{thm}

Theorem~\ref{det} is proved in Section~\ref{proofs}.
   
Consider a  pairing 
\be
s_\mu(z):\;
(L^q_{\La_1}\T\ldots\T L^q_{\La_n})_l\T
((L_{\La_1}\T\ldots\T L_{\La_n})_l)^*
\rightarrow \C,
\ee
defined by
\be
(f^{m_1}v^q_1\T\ldots\T f^{m_n}v^q_n)\T(f^{l_1}v_1\T\ldots\T f^{l_n}v_n)^*\,\mapsto \,
J_{\bar{l},\bar{m}}(z,\La)\in\C,
\ee
where $\bar{l},\bar{m}\in\Zb^n_l$ are $\La$-admissible.

\begin{corollary}\label{cor2}
Let $\Imag\mu\ne 0$.
Let $p\in\C$ and $z,\La\in\C^n$ satisfy conditions \Ref{step}-\Ref{weights3}. Then the
pairing
$s_\mu(z)$ is well defined and non-degenerate.
\end{corollary}

Corollary~\ref{cor2} follows from Theorem~\ref{det}.

\begin{corollary}\label{cor3}
Let $\Imag\mu\ne 0$.
Let  $p\in\C$ and $\La\in\C^n$ satisfy conditions \Ref{step}-\Ref{weights2} and
\Ref{step3}-\Ref{weights3}. Then any solution of the
qKZ equation with values in $(L_{\La_1}\T\ldots\T L_{\La_n})_l$ is a linear combination of
the hypergeometric solutions \Ref{Psi fin} with $p$-periodic coefficients.
\end{corollary}

Corollary~\ref{cor3} follows from Corollary~\ref{cor2}.

Similarly, for a permutation $\sigma\in{\Bbb S}^n$,
consider a  pairing
\be
s^\sigma_\mu(z):\;
(L^q_{\La_{\sigma_1}}\T\ldots\T L^q_{\La_{\sigma_n}})_l\T
((L_{\La_1}\T\ldots\T L_{\La_n})_l)^*
\rightarrow \C,
\ee
defined by
\be
(f^{m_{\sigma_1}}v^q_{\sigma_1}\T\ldots\T f^{m_{\sigma_n}}v^q_{\sigma_n})\T(f^{l_1}v_1\T\ldots\T f^{l_n}v_n)^*
\,\mapsto\,
c_{\bar{m}}(\La)I(w_{\bar{l}},W^\sigma_{\bar{m}})(z,\La)\in\C,
\ee
where $\bar{l},\bar{m}\in\Zb^n_l$ are $\La$-admissible.
   
\begin{corollary}\label{cor2 sigma}
Let $\Imag\mu\ne 0$.
Let $p\in\C$ and $z,\La\in\C^n$ satisfy conditions \Ref{step}-\Ref{weights3}. Then for any  permutation
$\sigma\in{\Bbb S}^n$, the pairing
$s^\sigma_\mu(z)$ is well defined and non-degenerate.
\end{corollary}

Corollary~\ref{cor2 sigma} follows from Corollaries~\ref{cor2} and \ref{transition functions}.

For a permutation $\sigma\in{\Bbb S}^n$ and
a $\La$-admissible index $\bar{m}\in\Zb^n_l$, consider a function
\bean\label{Psi fin sigma}
\Psi^\sigma_{\bar{m}}(z,\La)=\sum
I(w_{\bar{l}},W^\sigma_{\bar{m}})(z,\la) f^{l_1}v_1\T\ldots\T f^{l_n}v_n,
\eean
where the sum is over all $\La$-admissible $\bar{l}\in\Zb^n_l$.

\begin{corollary}\label{cor3 sigma}
Let $\Imag\mu\ne 0$. Let $\sigma\in{\Bbb S}^n$ be a permutation.
Let  $p\in\C$ and $\La\in\C^n$ satisfy conditions \Ref{step}-\Ref{weights2} and
\Ref{step3}-\Ref{weights3}. Then any solution of the
qKZ equation with values in $(L_{\La_1}\T\ldots\T L_{\La_n})_l$ is a linear combination of
the hypergeometric solutions \Ref{Psi fin sigma} with $p$-periodic coefficients.
\end{corollary}

Corollary~\ref{cor3 sigma} follows from Corollary~\ref{cor2 sigma}.

The map $s^\sigma_\mu(z)$ induces a map 
\bean\label{map sigma}
\tilde{s}^\sigma_\mu(z):\;
(L^q_{\La_{\sigma_1}}\T\ldots\T L^q_{\La_{\sigma_n}})_l
\rightarrow
(L_{\La_1}\T\ldots\T L_{\La_n})_l.
\eean
This map is an isomorphism of vector spaces for all $z$ satisfying \Ref{resonance}. For a fixed
vector $v^q\in (L^q_{\La_{\sigma_1}}\T\ldots\T L^q_{\La_{\sigma_n}})_l$, the
$(L_{\La_1}\T\ldots\T L_{\La_n})_l$-valued
function $\tilde{s}^\sigma_\mu(z)v^q$ is a solution of the qKZ equation. This construction
identifies  the space ${\cal 
S}$ of all meromorphic solutions to the qKZ equation with the space 
$(L^q_{\La_{\sigma_1}}\T\ldots\T L^q_{\La_{\sigma_n}})_l\T F$ where $F$ is the field of
scalar meromorphic functions in
$z_1,\dots,z_n$, $p$-periodic with respect to each of the variables,
\be
\iota_\sigma:\;
(L^q_{\La_{\sigma_1}}\T\ldots\T L^q_{\La_{\sigma_n}})_l\T F\to {\cal S}.
\ee

\bigskip

Let now $\Imag\mu=0$.
Consider a pairing
\be
s_\mu(z):\;
(L^q_{\La_1}\T\ldots\T L^q_{\La_n})^{sing}_l\T
((L_{\La_1}\T\ldots\T L_{\La_n})_l)^*
\rightarrow \C,
\ee
defined by
\be
v^q\T(f^{l_1}v_1\T\ldots\T f^{l_n}v_n)^*\,\mapsto \,
I(w_{\bar{l}},{\frak B}_{\rm id}(z,\La)v^q)(z,\La)\in\C,
\ee
where $\bar{l}\in\Zb^n_l$ is admissible, $v^q\in(L^q_{\La_1}\T\ldots\T L^q_{\La_n})^{sing}_l$. The map
${\frak B}_{\rm id}(z,\La)$ is defined in \Ref{trig.ident.1}, the integral
$I(w_{\bar{l}},{\frak B}_{\rm id}(z,\La)v^q)(z,\La)$ is taken over an integration contour associated to $(z,\La)$.

The map $s_\mu(z)$ induces a map
\bean\label{map'}
\tilde{s}_\mu(z):\;
(L^q_{\La_1}\T\ldots\T L^q_{\La_n})_l^{sing}
\rightarrow
(L_{\La_1}\T\ldots\T L_{\La_n})_l.
\eean

By Corollary~\ref{cor1'}, the image of $\tilde{s}_\mu(z)$ belongs to $(L_{\La_1}\T\ldots\T L_{\La_n})_l^{sing}$.

\begin{thm}\label{det'}
Let $\Imag\mu= 0$. 
Let $p\in\C$ and $z,\La\in\C^n$ satisfy conditions \Ref{step}-\Ref{weights3}. Let
$\sum\limits_{m=1}^n 2\La_m-s\not\in\{ p\Z_{<0}\}$ for $s=l-1,\dots, 2l-2 $.
Then the map
\be
\tilde{s}_\mu(z):\; 
(L^q_{\La_1}\T\ldots\T L^q_{\La_n})^{sing}_l
\rightarrow
(L_{\La_1}\T\ldots\T L_{\La_n})_l^{sing}
\ee
is an isomorphism of vector spaces.
\end{thm}

Theorem~\ref{det'} is proved similarly to Theorem~\ref{det} using Theorem 5.15 in \cite{TV1}.

\begin{corollary}\label{cor3'}
Let $\Imag\mu= 0$.
Let  $p\in\C$ and $\La\in\C^n$ satisfy conditions \Ref{step}-\Ref{weights2} and
\Ref{step3}-\Ref{weights3}.
Let $\sum\limits_{m=1}^n 2\La_m-s\not\in\{ p\Z_{<0}\}$ for $s=l-1,\dots, 2l-2 $.
Then any solution of
the qKZ equation with values in $(L_{\La_1}\T\ldots\T L_{\La_n})^{sing}_l$ is a linear
combination
of the hypergeometric solutions \Ref{Psi fin'} with $p$-periodic coefficients.  
\end{corollary}  

Corollary~\ref{cor3'} follows from Theorem~\ref{det'}.

\subsection{Asymptotic solutions}\label{as.sol.}

Let $U$ be a domain in $\C^n$ and let $M_1,\ldots,M_n$ be $sl_2$ modules. A basis $\Psi_1,\dots,\Psi_N$
of solutions to the qKZ equation
with values in $(M_1\T\ldots\T M_n)_l$, is
called an \emph{asymptotic solution} in the domain $U$ if 
\be
\Psi_j(z)=\exp(\sum_{m=1}^n a_{mj}z_m/p)
\prod_{1\le m<k\le n} (z_k-z_m)^{b_{jkm}}(v_j+o(1)),
\ee
where $a_{mj}$ and $b_{jkm}$ are suitable numbers, $v_1,\dots,v_N$ are vectors which form a basis
in $(M_1\T\ldots\T M_n)_l$ and $o(1)$ tends to zero as $z$ tends to infinity in $U$. The domain $U$
is called an \emph{asymptotic zone}.

Let $\Psi^1_1(z),\dots,\Psi^1_N(z)$ and $\Psi^2_1(z),\dots,\Psi^2_N(z)$ be asymptotic solutions
in asymptotic zones $U_1\in\C^n$ and $U_2\in\C^n$, respectively. Then
\be
\Psi^1_i (z)\,=\, \sum_{j=1}^N\, T^j_i(z)\, \Psi^2_j (z), \qquad i=1,..., N
\ee
for suitable meromorphic functions $T^j_i(z)$. The matrix valued function $T(z)=(T_i^j(z))$
is $p$ periodic with respect to variables $z_1,...,z_n$ and has nonzero determinant for
generic $z$. The function $T(z)$ is called  {\it
the transition function between the asymptotic solutions $\Psi^1$ and $\Psi^2$}.

We describe the asymptotic zones, asymptotic solutions and transition functions of the qKZ
equation with values in $(L_{\la_1}\T\ldots\T L_{\la_n})_l$. Our results are parallel to the
results of \cite{TV1}, where the asymptotic zones, asymptotic solutions and transition
functions are described for the  qKZ equation with values in $(V_{\la_1}\T\ldots\T
V_{\la_n})_l$.

For a permutation $\sigma\in{\Bbb S}^n$, define an asymptotic zone in $\C^n$,
\bean\label{zone}
U_\sigma=\{z\in\C^n\,|\,\Real z_{\sigma_1}\ll\ldots\ll\Real z_{\sigma_n}\}.
\eean
Say that $z\to\infty$ in $U_\sigma$ if $\Real (z_{\sigma_m}-z_{\sigma_{m+1}})\to -\infty$ for all
$m=1,\dots,n-1$.

For each permutation $\sigma\in{\Bbb S}^n$, let $\Psi^\sigma(\la)$ be the set of solutions $\{\Psi^\sigma_{\bar{m}}\}$
of the qKZ equation with
values in $(L_{\la_1}\T\ldots\T L_{\la_n})_l$ given by \Ref{Psi fin sigma}.
Here $\bar{m}\in\Zb^n_l$ runs through the set of all $\la$-admissible
indices. 

\begin{thm}\label{asympt}
Let $\Imag\mu\neq 0$.
Let $p\in\C,\La\in\C^n$ satisfy the conditions \Ref{step}-\Ref{weights2}.
Then for any permutation $\sigma\in{\Bbb S}^n$, the
set of solutions $\Psi^\sigma(\La)$ forms a basis of solutions. Moreover, the basis of
solutions $\Psi^\sigma(\La)$ is
an asymptotic solution
of the qKZ
equation with values in $(L_{\La_1}\T\ldots\T L_{\La_n})_l$ in the
asymptotic zone $U_\sigma$. Namely,
\bea
\lefteqn{\Psi^\sigma_{\bar{l}}(z)=\Xi_{\bar{l}}\,\exp(\mu\sum_{m=1}^nl_mz_m/p)\times}
\\&&
\times\prod_{1\le k<m\le n} ((z_{\sigma_k}-z_{\sigma_m})/p)^
{2(l_{\sigma_k}\La_{\sigma_m}+l_{\sigma_m}\La_{\sigma_k}-l_{\sigma_k}l_{\sigma_m})/p}
(f^{l_1}v_1\T\ldots\T f^{l_n}v_n+o(1)) 
\eea
as $z\to\infty$ in $U_\sigma$ so that at any moment condition \Ref{resonance} is satisfied.
Here $|\arg((z_{\sigma_k}-z_{\sigma m})/p)|<\pi$ and $\Xi_{\bar{l}}$ is a constant
independent of the permutation $\sigma$ and given by
\bea
\lefteqn{\Xi_{\bar{l}}=(2i)^l\, l!\, \Gamma(-1/p)^{-l}
\prod_{m=1}^n \left[ (e^\mu-1)^{(l_m(l_m-1)-2l_m\La_m)/p}\times\right.}
\\&&
\times\exp((\mu+\pi i)(l_m\La_m-l_m(l_m-1)/2)/p)\prod^{l_m-1}_{s=1}
\left.\Gamma((2\La_m-s)/p)\Gamma(-(s+1)/p)\right],
\eea
where  $0\le\arg(e^\mu-1)<2\pi$.
\end{thm}

The proof of Theorem~\ref{asympt} is the same as the proof of Theorem 6.4 in \cite{TV1}. 

Notice that the asymptotics of the basis of solutions $\Psi^\sigma(\La)$ determine the basis uniquely.
Namely, if a basis of solutions
meromorphically depends on parameters $\mu,z,\la$ and has asymptotics in $U_\sigma$ described in
Theorem~\ref{asympt}, then such a basis coincides with the basis $\Psi^\sigma(\La)$,
see the Remark after Theorem 6.4 in \cite{TV1}.

The isomorphism 
$\iota_\sigma:\;
(L^q_{\La_{\sigma_1}}\T\ldots\T L^q_{\La_{\sigma_n}})_l\T F\to {\cal S}$, defined in
Section~\ref{determinant},
allows us to identify the basis of solutions 
$\Psi^\sigma(\La)$ with a basis of $(L^q_{\La_{\sigma_1}}\T\ldots\T L^q_{\La_{\sigma_n}})_l\T F$,
considered as a vector
space over $F$. Then the transition function between two asymptotic solutions $\Psi^\sigma(\La)$ and
$\Psi^\nu(\La)$ is identified with an $F$-linear map
$(L^q_{\La_{\nu_1}}\T\ldots\T L^q_{\La_{\nu_n}})_l\T F\to
(L^q_{\La_{\sigma_1}}\T\ldots\T L^q_{\La_{\sigma_n}})_l\T F$. 

\begin{corollary}\label{trans}
Let $\Imag\mu\neq 0$.
Let $p\in\C,\La\in\C^n$ satisfy the conditions \Ref{step}-\Ref{resonance} and
\Ref{step3}-\Ref{weights3}. Then for any
permutation $\sigma\in{\Bbb S}^n$ and a simple transposition $(m,m+1)$, the transition function
between asymptotic solutions   $\Psi^\sigma(\La)$ and  $\Psi^{\sigma\circ(m,m+1)}(\La)$,
\be
(L^q_{\La_{\sigma_1}}\T\ldots\T L^q_{\La_{\sigma_{m+1}}}\T L^q_{\La_{\sigma_m}}\T\ldots\T L^q_{\La_{\sigma_n}})_l\T F
\to (L^q_{\La_{\sigma_1}}\T\ldots\T L^q_{\La_{\sigma_n}})_l\T F,
\ee
equals the operator
$P_{L_{\la_{\sigma_{m+1}}}^qL_{\la_{\sigma_m}}^q}R_{L_{\la_{\sigma_{m+1}}}^qL_{\la_{\sigma_m}}^q}
(e^{2\pi i(z_{\sigma_{m+1}}-z_{\sigma_m})/p})$ acting in the $m$-th and $(m+1)$-th factors.
\end{corollary}

Corollary~\ref{trans} follows from Theorem~\ref{asympt} and Corollary~\ref{transition functions}. 

Theorem~\ref{asympt} allows us to obtain another formula for the determinant of the hypergeometric pairing (cf. formula \Ref{det-prod}). 
\begin{thm}\label{det1} (Joint with V.Tarasov)
Let $\operatorname{Im}\,\mu\neq 0$. Let $p\in\C$ and $z,\La\in\C^n$ satisfy conditions \Ref{step}-\Ref{weights3}. Then the matrix $J_{\rm adm}^l(z,\La)$ defined in \Ref{adm.matrix} is non-degenerate. Moreover, for generic $p$,
\bea
\lefteqn{\det(J^l_{\rm adm}(z,\La))=}
\\&&
{\mathcal D}_l(\La)\exp(\mu\sum_{m=1}^n D_m(\La)z_m/p)\prod_{1\le k<m\le n}\,\prod_{s=0}^{d_{km}}\left(\frac{\Gamma((z_k+\la_k-z_m+\la_m-s)/p)}{\Gamma((z_k-\la_k-z_m-\la_m-s)/p)}\right)^{E_{km}(s,\La)},
\eea
where $d_{km}=\min\,\{l,\dim L_{\La_k}-1,\dim L_{\La_m}-1\}$, the functions ${\mathcal D}_l(\La),D_m(\La)$ and $E_{km}(s,\La)$ are given below.
\end{thm}

For $m=1,\dots,n$ and $\La\in\C^n$, set 
\be
D_m(\La)=\sum_{r=1}^{\min\{l,\dim L_{\La_m}-1\}}(r\dim (L_{\La_1}\T\ldots\T\widehat{L_{\La_m}}\T\ldots\T L_{\La_n})_{l-r}).
\ee

For $k,m=1,\dots,n$, $s\in\Z_{\ge 0}$ and $\La\in\C^n$, set
\be
E_{km}(s,\La)=\sum_r\,(\dim(L_{\La_k}\T L_{\La_m})_r-s-1)\,\dim(L_{\La_1}\T\ldots\T\widehat{L_{\La_k}}\T\ldots\T\widehat{L_{\La_m}}\T\ldots\T L_{\La_n})_{l-r},
\ee
where the sum is over $r=s+1,\dots,\min\,\{l,\dim L_{\La_k}+\dim L_{\La_k}-1-s\}$.

Finally, for $\La\in\C^n$, set 
\be
{\mathcal D}_l(\La)=\prod_{\bar{l}}\,c_{\bar{l}}\,\Xi_{\bar{l}},
\ee
where the product is over $\La$-admissible indices $\bar{l}\in{\mathcal Z}^n_l$.

The proof of Theorem~\ref{det1} is the same as the proof of Theorem 6.4 in \cite{TV1}.

\section{Proofs.}\label{proofs}
In this Section we collected the proofs of statements from Section~\ref{main}.

\subsection{Proof of Theorem~\ref{an.cont.1}}
Let $\Imag\mu\neq 0$.
Fix a rational weight function $w=w_{\bar{l}}$, and a trigonometric weight function $W=W_{\bar{m}}\in\G$,
such that $B(\bar{l})\bigcap B(\bar{m})=\emptyset$.

For $z,\la\in\C^n$ such that $\Real (z_i+\la_i)<0$ and
$\Real (z_i-\la_i)>0$ for all $i=1,\dots,n$, the integral $I(w,W)(z,\la)$ is defined by formula
\Ref{int},
\be
I(w,W)(z,\la)=\int\limits_{\Real t_i=0,\atop
i=1,\dots,l}f(t,z,\la)\,d^lt,
\ee
where $f(t,z,\la)=\Phi_l(t,z,\la)w(t,z,\la)W(t,z,\la)$.
The integrand $f$ has simple poles at 
\bean\label{poles1}
t_i=z_j\pm(\la_j+kp),\qquad k\in\Z_{\ge0},\qquad i=1,\dots,l,\qquad j=1,\dots,n,
\eean
and at 
\bean\label{poles2}
t_i-t_j=\pm(1-kp), \qquad k\in\Z_{\ge 0},\qquad i,j=1,\dots,l,\qquad i<j.
\eean
On the complex line the poles of the first type and the integration contour can be
represented as follows.

\setlength{\unitlength}{1cm}
\begin{picture}(14,4)
\thicklines
\put (7,0){\line(0,1){4}}
\put (7.1,0){$t_1,\dots,t_l$}  

\put (5.5,2.5){\circle*{0.1}}
\put (5.1,2.7){$z_1+\la_1$}
\put (3.5,2.5){\circle*{0.1}}
\put (2.7,2.7){$z_1+\la_1+p$}
\put (1.5,2.5){\circle*{0.1}}
\put (0.0,2.7){$z_1+\la_1+2p$}

\put (8.5,2.5){\circle*{0.1}}
\put (7.8,2.7){$z_1-\la_1$}
\put (10.5,2.5){\circle*{0.1}}
\put (9.6,2.7){$z_1-\la_1-p$}
\put (12.5,2.5){\circle*{0.1}}
\put (12.0,2.7){$z_1-\la_1-2p$}

\put (5.5,1){\circle*{0.1}}
\put (5.1,1.2){$z_2+\la_2$}
\put (3.5,1){\circle*{0.1}}
\put (2.7,1.2){$z_2+\la_2+p$}
\put (1.5,1){\circle*{0.1}}
\put (0.0,1.2){$z_2+\la_2+2p$}

\put (8.5,1){\circle*{0.1}}
\put (7.8,1.2){$z_2-\la_2$}
\put (10.5,1){\circle*{0.1}}
\put (9.6,1.2){$z_2-\la_2-p$}
\put (12.5,1){\circle*{0.1}}
\put (12.0,1.2){$z_2-\la_2-2p$}

\put (3,0){$\ldots$}
\put (10,0){$\ldots$}

\end{picture}
\medskip
\medskip

The integral $I(w,W)(z,\la)=I(w,W)(z,\la,p)$ is a meromorpic function of $z,\la,p$.

Fix $Z,\La\in\C^n$, $P\in\C$, satisfying conditions \Ref{step}-\Ref{resonance}. Our goal is to prove that
$I(w,W)(Z,\La,P)$ is well defined and is given by the integral over a suitable cycle.

We fix parameters $z^0,\la^0,p^0$ in such a way that $\Real\la^0_k<0$, 
$p^0$ is a negative number with large absolute value and
$|\Imag (z_k^0+\La_k-z_m^0-\La_m)|$ is large
for $k,m=1,\dots,n$ and $k\neq m$. Namely,
let $p^0=-2\sum\limits_{j=1}^n|\La_j|-1$, $z^0_k=i\,(-\Imag \La_k+3kA)$,
$\la^0_k=-1+i\,\Imag \La_k$, 
and $A$ is a large real number such that $A>2|\La_k|$, $k=1,\dots,n$.
 
The proof is done according to the following plan.
First, we represent the integral $I(w,W)(z^0,\la^0,p^0)$ as a sum of new integrals of different
dimensions. Then we analytically continue the function $I(w,W)(z,\la,p)$ from $z^0,\la^0,p^0$ to
$z^0,\La,p^0$. Then we analytically continue the function $I(w,W)(z,\La,p)$ from
$z^0,\La,p^0$ to $z^0,\La,P$. In the last step of the proof we analytically continue
the function
$I(w,W)(z,\La,P)$ from $z^0,\La,P$ to
$Z,\La,P$.

For $\bar{k}\in\Z_{\ge 0}^n$, set $k=k_1+\dots+k_n$ and define the multiple residues of $f$ by the
formula
\bean\label{mult.res}
\lefteqn{\res_{\bar{k}}f=\res_{t_k=z_n+\la_n-k_n+1}\ldots
\res_{t_{k_1+k_2}=z_2+\la_2-k_2+1}\ldots
\res_{t_{k_1+1}=z_2+\la_2}}\notag
\\ &&
\res_{t_{k_1}=z_1+\la_1-k_1+1}\ldots
\res_{t_2=z_1+\la_1-1}
\res_{t_1=z_1+\la_1}\,f.
\eean

If $k>l$, then we set $\res_{\bar{k}}f=0$.

The function $\res_{\bar{k}}f$ is a function of variables $t_{k+1},\dots,t_l;z,\la$.
In particular, for $\bar{k}=(0,\dots,0)$, we have $\res_{\bar{k}} f=f$.

\begin{lemma}\label{order}
Let $\bar{k}, \bar{k}^\prime\in\Z^n_{\ge0}$. If $\bar{k}^\prime$ can be obtained by a permutation of
components of $\bar{k}$
then $\res_{\bar{k}}f=\res_{\bar{k}^\prime}f$.
\end{lemma}

Lemma~\ref{order} follows from the symmetry of the integrand $f$, see Remark in
Section~\ref{actions}.

\begin{lemma}\label{double res}
For all $i,j=1,\dots,l$, and for any $a\in \C$, we have
\be
\res_{t_j=a}\,\res_{t_i=a}\,f=0.
\ee
\end{lemma}
\begin{proof}
Recall that $f=\Phi_l wW$.
On the hyperplane $t_i=t_j$, any  function which is symmetric with respect to either
rational or trigonometric action is identically equal to zero. In particular, $w$ is
equal to zero at
the hyperplane $t_i=t_j$. Since all poles of $f$ are
at most of first order, we have the statement of the Lemma.
\end{proof}

\begin{lemma}\label{res}
Let $ \bar{k}\in\Z^n_{\ge0}$.
The function $\res_{\bar{k}}f$ is
holomorpic for all $t_{k+1},\dots,t_l;z,\la,p$ satisfying
conditions \Ref{step} and such that
\be
t_i\neq z_j\pm(\la_j+sp),\qquad s\in\Z_{\ge0},\qquad i=k+1,\dots,l,\qquad j=1,\dots,n,
\ee
\be
t_i-t_j\neq\pm(1-sp), \qquad s\in\Z_{\ge 0},\qquad i,j= k+1,\dots,l,\qquad i<j,
\ee
\be
t_i\neq z_j+\la_j-m\pm(1-sp),
\;\; s\in\Z_{\ge 0},\;\; i=k+1,\dots,l,\;\; j= 1,\dots,n,\;\; m=0,\ldots,k_j,
\ee
\be
z_i+\la_i-z_j\pm(\la_j+sp)\neq m,
\qquad s\in\Z_{\ge 0},\qquad i,j=1,\dots,k,\qquad m=0,\ldots,k,\qquad i\neq j,
\ee
\be
2\la_i+sp\neq m,
\qquad s\in\Z_{\ge 0},\qquad i=1,\dots,k,\qquad m=0,\ldots,k_i.
\ee
\end{lemma}

Lemma~\ref{res} follows from formulas \Ref{poles1},\Ref{poles2} for the poles of $f$,
and formula \Ref{mult.res} for
the function $\res_{\bar{k}}f$.

For $u\in\R^n$, define a curve  ${\cal C}_u=\{{\cal C}_u(x)\in\C\,|\,x\in\R\}$ with
the following properties. The
curve ${\cal C}_u$ consists of $2n+1$ line segments and none of the line segments is horizontal (i.e. 
no line-segment is parallel to the
real line). There are two line segments (rays) which go to $+i\infty$ and $-i\infty$,
these two rays are parts of the imaginary
line. There are $n$ vertical segments with real coordinates $u_1,\dots,u_n$ and
imaginary
coordinates close to $\Imag (z^0_1+\la^0_1),\dots, \Imag (z_n^0+\la_n^0)$. 

The precise formulas defining ${\cal C}_u$  are as follows.
For $x\le A$
and for
$x\ge(3n-1)A$, let
${\cal C}_u(x)=ix$.
For 
$(3m-1)A\le x\le (3m+1)A$, $m=1,\dots,n$, let ${\cal C}_u(x)=u_m+ix$, and this segment is called the
\emph{$m$-th
vertical segment}. For
$m=0,\dots,n$, let
$\{{\cal C}_u(x), (3m+1)A\le x\le (3m+2)A\}$ be a parametrization of the line segment connecting ${\cal
C}_u((3m+1)A)$
and ${\cal C}((3m+2)A)$.

See an example of such a curve on the Figure of the proof of Lemma~\ref{la-cont}.

\begin{lemma}\label{contour}
Let $,\La,z^0,\la^0,p^0$ be as above. For any sufficiently small $\epe$,
\bean\label{sum}
I(w,W)(z^0,\la^0,p^0)=\sum_{\bar{k}}\,\frac{l!}{n!\,k_1!\dots
k_n!(l-k)!}\int\limits_{t_j\in{\cal C},\atop
j=k+1,\dots,l}\res_{\bar{k}}f(t,z^0,\la^0,p^0)\,d^{l-k}t,
\eean
where $d^{l-k}t=dt_{k+1}\ldots dt_l$. Here ${\cal C}$ denotes the curve ${\cal C}_u$
with
$u_m=\min\{-\Real(2\La_m-\la_m^0)-\epe,0\}$, $m=1,\dots,n$,
and the sum is over all $\bar{k}\in\Z^n_{\ge 0}$ such that
$k\le l$, and $k_m\le 2\Real\La_m+1$ for $m=1,\dots,n$.
\end{lemma}

\begin{proof}
The integrand decays exponentially with respect to all $t_1,\dots,t_l$. We move the contour of
integration with respect to each variable
$t_i$ to the left avoiding the poles at $t_i-t_j=(1-kp)$. 

Initially, each of the variables runs through the imaginary line. We deform the contours of
integration with respect to $t_1,\dots,t_l$ in such a way that at every moment the contours of integration with respect to
$t_1,\dots,t_l$ are curves ${\cal C}_{u^i}$ with
different parameters $u^i\in\R^n$.

First, we separate the contours of integration with respect to different variables by
small distances as follows.

\be
I(w,W)(z^0,\la^0,p^0)=
\int\limits_{
\begin{array}{c}
\scriptstyle{\Real t_j=0,}
\\
\scriptstyle{j=i,\dots,l}
\end{array}}
\!\!\!fd^lt=
\int\limits_{
\begin{array}{c}
\scriptstyle{t_j\in{\cal C}_{u^j},}\\
\scriptstyle{j=1,\dots,l}
\end{array}}
\!\!\!fd^lt,
\ee
where $u^j=((j-l-1)\rho,\dots,(j-l-1)\rho)\in\R^n$, and $\rho$ is a small positive
number. 

We deform the integration cycle $\{t\in\C^l\,|\,t_j\in{\cal C}_{u^j},\, j=1,\dots,l\}$ changing
$u^j,\, j=1,\dots,l$. Namely,
we move the vertical segments of contours of integrations ${\cal C}_{u^j}$,
$j=1,\dots,n$, to the left,
until we reach the vertical segments of the curve $\cal C$, described in the Lemma,
and get the integration cycle $\{t\in\C^l\,|\,t_j\in{\cal C},\, j=1,\dots,l\}$. During the deformation, we
keep the same small distances between the
contours of integrations ${\cal C}_{u^j}$ with respect to different variables. At every moment, when one of
the vertical segments of the contour of integration
goes through a pole of the integrand, we add another integral of lower dimension which
is the integral
with respect to remaining variables of
the residue. For example, the first such event leads to the following decomposition.

\bea
\lefteqn{\int\limits_{
\begin{array}{c}
\scriptstyle{t_j\in{\cal C}_{u^j},}\\
\scriptstyle{j=1,\dots,l}
\end{array}}
\!\!\!fd^lt=
\int\limits_{
\begin{array}{c}
\scriptstyle{t_j\in{\cal C}_{(\Real \la_1^0+j\rho,u^j_2,\dots,u^j_n)},}\\
\scriptstyle{j=1,\dots,l}
\end{array}}
\!\!\!fd^lt=}
\\&&
\int\limits_{
\begin{array}{c}
\scriptstyle{t_1\in{\cal C}_{(\Real \la_1^0-\rho/2},u^1_2,\dots,u^1_n)}\\
\scriptstyle{t_j\in{\cal C}_{(\Real \la_1^0+(j-3/2)\rho,u^j_2,\dots,u^j_n)},}\\
\scriptstyle{j=2,\dots,l}
\end{array}}
\!\!\!fd^lt+
\int\limits_{
\begin{array}{c}
\scriptstyle{t_j\in{\cal C}_{(\Real \la_1^0+(j-3/2)\rho,u^j_2,\dots,u^j_n)},}\\  
\scriptstyle{j=2,\dots,l}
\end{array}}
\!\!\!\res_{t_1=z_1+\la_1}fd^{l-1}t.
\eea

Recall that $\Real \la_1^0=-1$.
In the first term of the resulting decomposition, moving the first vertical segment of
the contour ${\cal C}_{(\Real \la_1^0+1/2\rho,u^2_2,\dots,u^2_n)}$ corresponding
to $t_2$,
through the point $z^0_1+\la^0_1$, we also get a lower dimensional integral. This integral is equal to
the
second term of the above decomposition by Lemma~\ref{order}.
Note that due to Lemma~\ref{double res}, the integrand of the second term of the
decomposition does not have poles at
$t_j=z_1+\la_1$ for all
$j=2,\dots,l$. Thus, we can move the first vertical segments of the contours corresponding
to $t_j$, $j=2,\dots,l$,
through $z^0_1+\la^0_1$ without getting new residues.

Finally, we have
\bea
\lefteqn{I(w,W)(z^0,\la^0,p^0)=}
\\&&
\int\limits_{
\begin{array}{c}
\scriptstyle{t_j\in{\cal C}_{(\Real \la_1^0+(j-l-1)\rho,u^j_2,\dots,u^j_n)},}\\
\scriptstyle{j=1,\dots,l}
\end{array}}   
\!\!\!fd^lt+
l\int\limits_{
\begin{array}{c}
\scriptstyle{t_j\in{\cal C}_{(\Real\la_1^0+(j-l-1)\rho,u^j_2,\dots,u^j_n)},}\\   
\scriptstyle{j=2,\dots,l}  
\end{array}}
\!\!\!\res_{t_1=z_1+\la_1}fd^{l-1}t.
\eea

In the first summand we move the
contour of integration to 
\be
\{t\in\C^l\,|\, t_j\in{\cal
C}_{(\min\{-\Real(2\La_1-\la_1^0)-\epe,0\},u^j_2,\dots,u^j_n)},\,
j=1,\dots,l\}
\ee
and  the contour does not meet any other poles of the integrand.

The integrand of the second summand has 
a simple pole at  $t_j=z_1+\la_1-1$, since $f$ had a pole at 
$t_1-t_j= 1$. It also may have a pole at $2\la_1=0$ coming from the pole at 
$t_1=z_1-\la_1$, see Lemma~\ref{res}.
So, as before, we have
\be
\int\limits_{
\begin{array}{c}
\scriptstyle{t_j\in{\cal C}_{(\Real \la_1^0+(j-l-1)\rho,u^j_2,\dots,u^j_n)},}\\
\scriptstyle{j=2,\dots,l}
\end{array}}
\!\!\!\res_{t_1=z_1+\la_1}fd^{l-1}t=
\int\limits_{
\begin{array}{c}
\scriptstyle{t_j\in{\cal C}_{(\Real \la_1^0-1+(j-1)\rho,u^j_2,\dots,u^j_n)},}\\
\scriptstyle{j=2,\dots,l}
\end{array}}
\!\!\!\res_{t_1=z_1+\la_1}fd^{l-1}t=
\ee
\bea
\lefteqn{
\int\limits_{
\begin{array}{c}
\scriptstyle{t_j\in{\cal C}_{(\Real \la_1^0-1+(j-l-1)\rho,u^j_2,\dots,u^j_n)},}\\
\scriptstyle{j=2,\dots,l}
\end{array}}
\!\!\!\res_{t_1=z_1+\la_1}fd^{l-1}t+}
\\&& 
+(l-1)\int\limits_{
\begin{array}{c}
\scriptstyle{t_j\in{\cal C}_{(\Real \la_1^0-1+(j-l-1))\rho,u^j_2,\dots,u^j_n)},}\\
\scriptstyle{j=3,\dots,l}
\end{array}}
\!\!\!\res_{t_2=z_1+\la_1-1}\res_{t_1=z_1+\la_1}fd^{l-2}t.
\eea

Again, we move the contour of integration in the first summand to
\be
\{(t_2,\dots,t_l)\in\C^{l-1}\,|\, t_j\in{\cal
C}_{(\min\{-\Real(2\La_1-\la_1^0)-\epe,0\},u^j_2,\dots,u^j_n)},
\,j=2,\dots,l\}
\ee
without meeting any poles of the integrand,
and the integrand of the 
second summand has simple poles at  $t_j=z_1+\la_1-2$ coming from the poles at  
$t_2-t_j=1$. It may also have poles at $2\la_1=1$ coming from $t_2=z_1-\la_1$.

Continuing the process we prove Lemma~\ref{contour}.
\end{proof}

\begin{lemma}\label{la-cont}
The analytic continuation of the function $I(w,W)(z,\la,p)$ from $z^0,\la^0,p^0$ to a
neighborhood of the point $z^0,\La,p^0$ is given by
\bean\label{La-sum}
I(w,W)(z^0,\la,p^0)=\sum_{\bar{k}}\,\frac{l!}{n!\,k_1!\dots
k_n!(l-k)!}\int\limits_{t_i\in{\cal C},\atop
i=k+1,\dots,l}\res_{\bar{k}}f(t,z^0,\la,p^0)\,d^{l-k}t,
\eean
where $k=k_1+\ldots+k_n$, $d^{l-k}t=dt_{k+1}\ldots dt_l$. Here ${\cal C}$ denotes the curve ${\cal C}_u$
with
$u_m=\min\{-\Real\La_m-\epsilon,0\}$, $m=1,\dots,n$, $\epe$ is a sufficiently small positive number,
and the sum is over all $\bar{k}\in\Z^n_{\ge 0}$ such that
$k\le l$, and $k_m\le 2\Real\La_m+1$ for $m=1,\dots,n$. Morever if $B(\bar{l})\bigcap B(\bar{m})=\emptyset$, then sum \Ref{La-sum} is
well defined at $z^0,\La,p^0$.
\end{lemma}
\begin{proof}
We analytically continue each summand in \Ref{sum} from $z^0,\la^0,p^0$ to $z^0,\La,p^0$. 
By Lemma~\ref{contour}, the poles of the integrand of a summand and the contour of integration can be
pictured as follows.

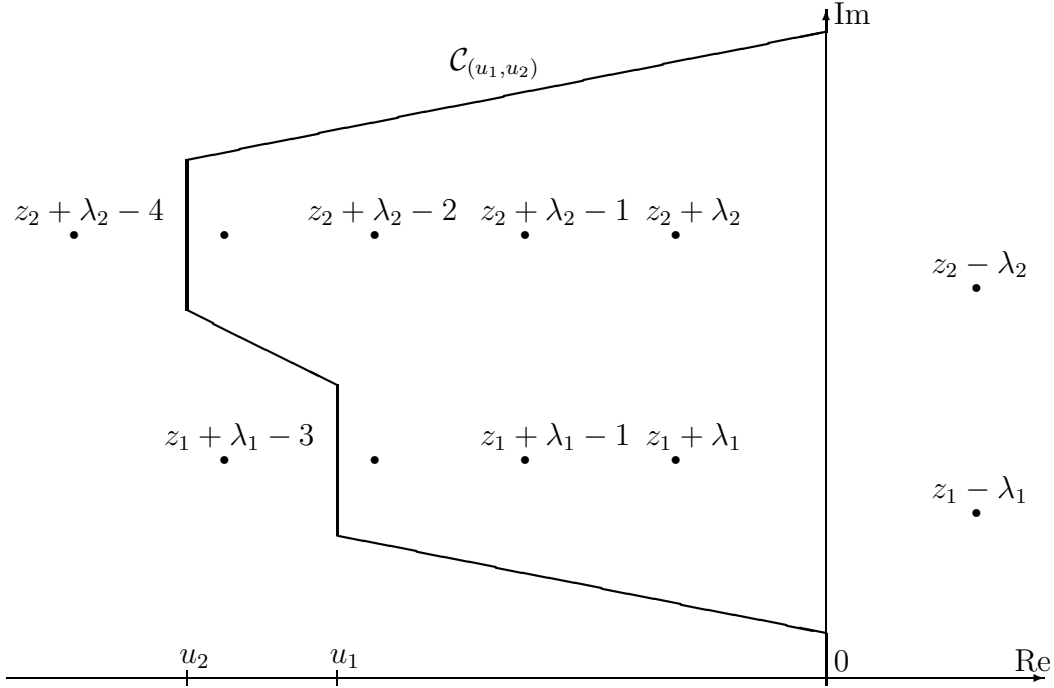
\setlength{\unitlength}{1cm}
\begin{picture}(14,9)

\thicklines
\put (11,0){\line(0,1){0.7}}
\put (11,8.7){\line(0,1){0.3}}
\put (4.5,2){\line(5,-1){6.5}}
\put (4.5,2){\line(0,1){2}}
\put (2.5,5){\line(0,1){2}}
\put (2.5,5){\line(2,-1){2}}
\put (2.5,7){\line(5,1){8.5}}

\thinlines
\put (11,0){\vector(0,1){9}}
\put (0.1,0.1){\vector(1,0){13.8}}

\put (4.5,0.0){\line(0,1){0.2}}
\put (2.5,0.0){\line(0,1){0.2}}

\put (11.1,8.8){Im}
\put (13.5,0.2){Re}
\put (2.4,0.3){$u_2$}
\put (4.4,0.3){$u_1$}
\put (11.1,0.2){0}
\put (6,8.2){${\cal C}_{(u_1,u_2)}$}

\put (13,2.3){\circle*{0.1}}
\put (12.4,2.5){$z_1-\lambda_1$}

\put (13,5.3){\circle*{0.1}}  
\put (12.4,5.5){$z_2-\lambda_2$}

\put (1,6){\circle*{0.1}}   
\put (3,6){\circle*{0.1}}
\put (5,6){\circle*{0.1}}
\put (7,6){\circle*{0.1}}
\put (9,6){\circle*{0.1}}   

\put (3,3){\circle*{0.1}}
\put (5,3){\circle*{0.1}}
\put (7,3){\circle*{0.1}}
\put (9,3){\circle*{0.1}}

\put (8.6,6.2){$z_2+\lambda_2$}
\put (6.4,6.2){$z_2+\lambda_2-1$}
\put (4.1,6.2){$z_2+\lambda_2-2$}
\put (0.2,6.2){$z_2+\lambda_2-4$}   

\put (8.6,3.2){$z_1+\lambda_1$}
\put (6.4,3.2){$z_1+\lambda_1-1$}
\put (2.2,3.2){$z_1+\lambda_1-3$}

\end{picture}

Here $z_j=z_j^0$, $\la_j=\la_j^0$, $j=1,\dots,n$.  

Note that there are also poles of the integrand of the type $t_i=z_j^0\pm(\la_j+p^0s)$, $s=1,2,\dots\,.$
These poles are far away from our picture according to our choice of $p^0$. Note also that
all poles are simple.

We move the parameters
$\la_m$ from $\la_m^0$ to $\La_m$. The contour of integration of each summand in Lemma~\ref{contour},
$\{t\in\C^l\,|\,t_j\in{\cal C}={\cal C}_{u(\la)},\,j=k+1,\dots,l\}$, depends on $\la$.
At every moment of deformation of parameters
$\la_m$, we define $I(w,W)(z^0,\la,p^0)$ by
\be
I(w,W)(z^0,\la,p^0)=\sum_{\bar{k}}\,\frac{l!}{n!\,k_1!\dots
k_n!(l-k)!}\int\limits_{t_i\in{\cal C}_{u(\la)},\atop
i=k+1,\dots,l}\res_{\bar{k}}f(t,z^0,\la,p^0)\,d^{l-k}t,
\ee
where $u(\la)=(u_1(\la_1),\dots,u_n(\la_n))\in\R^n$,
$u_m(\la_m)=\min\{-\Real(2\La_m-\la_m)-\epe,0\}$, $m=1,\dots,n$,
and the sum is over all $\bar{k}\in\Z^n_{\ge 0}$ such that
$k\le l$, and $k_m\le 2\Real\La_m$ for $m=1,\dots,n$.

If $\La_m\not\in\La^+$ for all $m=1,\dots,n$, then at
every moment of the deformation the integrand is holomorphic by Lemma~\ref{res}. On the
drawing, the points $z_i-\la_i$,
$i=1,\dots,n$, move to the left and the rest of the picture moves to the right. Note
that the points $z_i-\la_i$, $i=1,\dots,n$, encounter neither other points on the
picture nor the curve of integration ${\cal C}_{u(\la)}$. 
 
If, for some $i\in\{1,\dots,n\}$, $\La_i\in\La^+$, then by Lemma~\ref{res}, the function
$\res_{\bar{k}}(z,\la,p)$ can
have a pole at $\la_i=\La_i$ for $\bar{k}$ such that $k_i=2\La_i$. In this case, in the 
picture, in the last moment of the deformation the
point $z_i-\la_i$ coincides with the point $z_i+\la_i-k$.

However, if
$w=w_{\bar{l}},\; W=W_{\bar{m}}$ and $B(\bar{l})\bigcap B(\bar{m})=\emptyset$, then the pole is trivial by the
following two lemmas.

\begin{lemma}\label{main residue}
Let $i\in\{1,\dots,n\}$ and $\La_i\in\La^+$.
Let $W=W_{\bar{m}}$ and $m_i\le 2\La_i$. Then for any $w=w_{\bar{l}}$,
\be
\res_{\la_i=\La_i}
\res_{t_{2\La_i+1}=z_i+\la_i-2\La_i}\ldots\res_{t_2=z_i+\la_i-1}\res_{t_1=z_i+\la_i} 
\Phi_lwW=0.
\ee
\end{lemma}
\begin{proof}
Let i=1. We prove that
\be 
\res_{t_{2\La_1+1}=z_1+\la_1-2\La_1}\ldots\res_{t_2=z_1+\la_1-1}\res_{t_1=z_1+\la_1}
\Phi_lwW=0.
\ee
Recall that $W$ has the form
$W=\Sum_{\sigma\in {\Bbb S}^l}[\ldots]_\sigma^{trig}$. 
The residue is nontrivial only for the terms corresponding to permutations $\sigma$,
which satisfy 
the condition: \be
m_1\ge(\sigma^{-1})_1>(\sigma^{-1})_2>\ldots>(\sigma^{-1})_{2\La+1}.
\ee
If $m_1\le 2\La$ then there are no such permutations.

The cases $i=2,\dots,n$ are proved similarly.
\end{proof}

\begin{lemma}
Let $i\in\{1,\dots,n\}$ and $\La_i\in\La^+$.
Let $w=w_{\bar{l}}$ and $l_i\le 2\La_i$. Then for any $W=W_{\bar{m}}$,
\be
\res_{\la_i=\La_i}
\res_{t_{2\La_i+1}=z_i+\la_i-2\La_i}\ldots\res_{t_2=z_i+\la_i-1}\res_{t_1=z_i+\la_i}
\Phi_lw_{\bar{l}}W_{\bar{m}}=0.
\ee
\end{lemma}
\begin{proof}
Let i=1. We prove that
\be
\res_{t_{2\La_1+1}=z_1+\la_1-2\La_1}\ldots\res_{t_2=z_1+\la_1-1}\res_{t_1=z_1+\la_1}
\Phi_lw_{\bar{l}}\tilde{W_{\bar{m}}}=0.
\ee
Recall that $w_{\bar{l}}$ has the form
$w_{\bar{l}}=\Sum_{\sigma\in {\Bbb S}^l}[\ldots]_\sigma^{rat}$.
The resudue is nontrivial only for the terms corresponding to permutations $\sigma$
which satisfy the condition: 
\be
l_1\ge(\sigma^{-1})_1>(\sigma^{-1})_2>\ldots>(\sigma^{-1})_{2\La+1}.
\ee
Since $1\not \in B(\bar{l})$, we have $l_1\le 2\La$ and there are no such permutations.

The cases $i=2,\dots,n$ are proved similarly.
\end{proof}

Lemma~\ref{la-cont} is proved.
\end{proof}

Note that a residue of a meromorphic function can be described as an integral over a small circle.

For $B\in\C$ and $\epe>0$, denote ${\cal D}_B\subset\C$ the circle with center $B$ and radius $\epe$.
Denote ${\cal D}$ the curve ${\cal C}_u$ with
$u_m=\min\{-\Real\La_m-\epsilon,0\}$.

\begin{lemma}\label{p-cont}   
The analytic continuation of the function $I(w,W)$ from $z^0,\La,p^0$ to $z^0,\La,P$ is well defined. Moreover,
for sufficiently small $\epe>0$,
\bean\label{p-sum}
\lefteqn
{I(w,W)(z^0,\La,P)=}
\\&&
\sum_{m=0}^l\sum_{a_1,\dots,a_m}\,\frac{l!}{(l-m)!}
\int\limits_{t_i\in{\cal D},\atop i=m+1,\dots,l}
\left(\sum_{\sigma\in{\Bbb S}^m}\int\limits_{t_i\in{\cal D}_{a_i},\atop
i=1,\dots,m}f(t,z^0,\La,P)\,dt_{\sigma_1}\dots dt_{\sigma_m}\right)dt_{m+1}\dots dt_l,\notag
\eean
where the second sum is over all $a=(a_1,\dots,a_m)\in\C^m$ such that
for each $i=1,\dots,m$, there exist 
$j\in\{1,\dots,n\}$, $k,s\in\Z_{>0}$, such that $a_i=z_j+\La_j-k+sP$; $a_i$ is located to the right from ${\cal
D}$ and $a_i\neq
a_k$ for all $k,i=1,\dots,m$, $k\neq i$.	
\end{lemma}

Note that formula \Ref{p-sum} is a generalization of formula \Ref{La-sum}, however, in sum
\Ref{p-sum} there are many additional zero terms.

\begin{proof}
We analytically continue each summand in \Ref{La-sum}  with respect to $p$ from
$p^0$ to $P$ preserving $z^0,\La$. 
We move $p$ and
preserve the contour of integration $\cal I$ as long as the integrand is holomorphic for all
$t\in{\cal I}$ at $z,\La,p$ (cf. Lemma~\ref{res}). When a pole of the integrand goes through the contour of integration,
we do the same
procedure as in the proof of Lemma~\ref{contour}.

An individual integral on the right hand side of \Ref{p-sum} is over the cycle
\be
\{t\,|\,t_i\in{\cal D}_{a_i},\,i=1,\dots,m;\,t_i\in{\cal D},\,i=m+1,\dots,l\}.
\ee
The points $a_i$ are the poles of the integrand which were on the left of ${\cal D}$ at
the beginning of the deformation and are on the right of ${\cal D}$ at the final moment
of the deformation. All poles of the integrand which were on the right of ${\cal D}$ at
the beginning of the deformation, remain on the right at the every moment of the
deformation. The new poles $a_i$ on the right of ${\cal D}$ do not coincide with the
old poles on the right of ${\cal D}$ due to conditions \Ref{step}-\Ref{weights2} on
$P$ and $\La$.

Note, that the procedure of analytic continuation, described in Lemma~\ref{contour}, gives
decomposition \Ref{p-sum} in which the sum is over $a=(a_1,\dots,a_m)$, where some of
coordinates $a_i$ could be equal. However, such integrals are equal to zero. In fact,
let $a_i=a_j$ for some $i<j$. When we integrate over the variables $t_1,\dots,t_{j-1}$,
the resulting integrand could have a pole at $t_j=a_j$ of at most first order. Then the
reason of Lemma~\ref{double res} shows that the resulting integrand is holomorphic at
$t_j=a_j$.
\end{proof}

Let, ${\cal D}^\prime$ be the line $\{x\in\C\,|\,\Real x=\min\{\Real (Z_i-\La_i), \,
i=1,\dots,n\}-\epe\}$.

\begin{lemma}\label{z-cont}
The analytic continuation from $z^0,\La,P$ to $Z,\La,P$ is well defined. Moreover, for
sufficiently small $\epe>0$,
\bean\label{z-sum}
\lefteqn
{I(w,W)(Z,\La,P)=}
\\&&
\sum_{m=0}^l\sum_{b_1,\dots,b_m}\,\frac{l!}{(l-m)!}
\int\limits_{t_i\in{\cal D}^\prime,\atop i=m+1,\dots,l}
\left(\sum_{\sigma\in{\Bbb S}^m}\int\limits_{t_i\in{\cal D}_{b_i},\atop
i=1,\dots,m}f(t,Z,\La,P)\,dt_{\sigma_1}\ldots dt_{\sigma_m}\right)dt_{m+1}\ldots,dt_l,\notag
\eean
where the second sum is over all $b=(b_1,\dots,b_m)\in\C^m$ such that
for each $i=1,\dots,m$, there exist
$j\in\{1,\dots,n\}$, $k,s\in\Z_{>0}$, such that $b_i=Z_j+\La_j-k+sP$; and $b_i$ is located to the right from
${\cal D}^\prime$.
\end{lemma}

\begin{proof}

First, we move the parameters
$z_m$ from $z_m^0$ to $z^0_m+\Real Z_m$. The contour of integration of a summand in \Ref{p-sum} has the
form
\be
\{t\in\C^l\,|\,t_i\in{\cal D}_{a_i(z)},\,i=1,\dots,m;\,t_i\in{\cal D}={\cal C}_{u(z)},\,i=m+1,\dots,l\},
\ee
where $u(z)=\Real (z_m-\La_m)-\epe$, and depends on $z$. At every
moment of the deformation of parameters $z_m$, we define $I(w,W)(z,\La,P)$ by
\bea
\lefteqn
{I(w,W)(z,\La,P)=}
\\&&
\sum_{m=0}^l\sum_{a_1,\dots,a_m}\,\frac{l!}{(l-m)!}
\int\limits_{t_i\in{\cal C}_{u(z)},\atop i=m+1,\dots,l}
\left(\sum_{\sigma\in{\Bbb S}^m}\int\limits_{t_i\in{\cal D}_{a_i},\atop
i=1,\dots,m}f(t,z,\La,P)\,dt_{\sigma_1}\ldots
dt_{\sigma_m}\right)dt_{m+1}\ldots,dt_l,\notag
\eea
where the second sum is over all $a=(a_1,\dots,a_m)\in\C^m$ such that
for each $i=1,\dots,m$, there exist
$j\in\{1,\dots,n\}$, $k,s\in\Z_{>0}$, such that $a_i=z_j+\La_j-k+sP$; and $a_i$
is located to the right from ${\cal C}_{u(z)}$.

Consider the summand related to $a\in\C^m$. We move the contour of integration with
respect to $t_{m+1},\dots,t_l$ from ${\cal C}_{u(z^0+\Real Z)}$ to ${\cal D}^\prime$.
We use the method described in Lemma~\ref{contour}.

Finally, we move parameters $z_1\dots,z_n$ from $z_1^0+\Real Z_1,\dots,z_n^0+\Real Z_n$ to
$Z_1\dots,Z_n$ in the same way as we moved $\la_1,\dots,\la_n$ in
Lemma~\ref{la-cont}.

An individual integral in the RHS of \Ref{z-sum} is over the cycle
\be
\{t\,|\,t_i\in{\cal D}_{b_i},\,i=1,\dots,m;\,t_i\in{\cal D}^\prime,\,i=m+1,\dots,l\}.
\ee
The points $b_i$ are the poles of the integrand which were on the left of ${\cal D}$ at
the beginning of the deformation and are on the right of ${\cal D}^\prime$ at the final
moment of the deformation. All poles of the integrand which were on the right of
${\cal D}$ at
the beginning of the deformation, remain on the right of ${\cal D}^\prime$ at the every
moment of the deformation. The new poles $b_i$ on the right of ${\cal D}^\prime$ do not
coincide with the old poles on the right of ${\cal D}^\prime$ due to conditions
\Ref{resonance} on $Z$.
\end{proof}
Theorem~\ref{an.cont.1} is proved.
$\;\Box$

\subsection{Proof of Theorems~\ref{an.cont.2}, \ref{functor}}

Let $\Imag\mu\neq 0$.
Fix a rational weight function $w=w_{\bar{l}}$, and a trigonometric weight function $W=W_{\bar{m}}\in\G$.
Fix $\La\in\C^n$.

We fix parameters $z,p,\la^0$ in such a way that $\Real \la^0_m<0$,
$p$ is a negative number with large absolute value and
$|\Imag (z_k+\La_k-z_m-\La_m)|$ is large
for $k,m=1,\dots,n$, $k\neq m$. Namely,
let $p=-2\sum\limits_{j=1}^n|\La_j|-1$, $\la_k^0=-1+i\Imag \La_k$, $z^0_k=i\,(-\Imag \La_k+3kA)$,
where $A$ is a large real number such that $A>2|\La_k|$, $k=1,\dots,n$.
   
By Lemma~\ref{la-cont},
\bean\label{la-sum1}
I(w,W)(z,\la,p)=\sum_{\bar{k}}\,\frac{l!}{n!\,k_1!\dots
k_n!(l-k)!}\int\limits_{t_i\in{\cal C},\atop
i=k+1,\dots,l}\res_{\bar{k}}f(t,z,\la,p)\,d^{l-k}t,
\eean
where $k=k_1+\ldots+k_n$, ${\cal C}$ denotes the curve ${\cal C}_u$ with
$u_m=\min\{-\Real(2\La_m-\la_m)-\epsilon,0\}$, $m=1,\dots,n$, $\epe$ is a sufficiently small positive number,
and the sum is over all $\bar{k}\in\Z^n_{\ge 0}$ such that
$k\le l$, and $k_m\le 2\Real\La_m+1$ for $m=1,\dots,n$.

Let $i\in\{1,\dots,n\}$. We have the following cases.

1). If $i$ does not belong to $B_{\La}(\bar{m})$, then the analytic continuation of the function $I(w,W)(z,\la,p)$
with respect to $\la_i$ from
$\la_i^0$ to $\La_i$ is well defined, see the proof of Lemma~\ref{la-cont}. 

2). If $i$ belongs to $B_{\La}(\bar{m})$ and $i$ does not belong to $B_{\La}(\bar{l})$, then the analytic
continuation of the function
$I(w,W)(z,\la,p)$ with respect to $\la_i$ from $\la_i^0$ to $\La_i$ is well defined, see the proof of
Lemma~\ref{la-cont}. In this case 
$c_{\bar{m}}(\la,p)$ is equal to zero at $\la_i=\La_i$, so $J_{\bar{l},\bar{m}}(z,\la,p)=c_{\bar{m}}(\la,p)I(w,W)(z,\la,p)$ is equal to
zero at $\la_i=\La_i$.

3). If $i$ belongs to $B_{\La}(\bar{m})$ and to $B_{\La}(\bar{l})$, then in decomposition
\Ref{la-sum1} all the summands are well defined at $\la_i=\La_i$ except for the summands
corresponding to
$\bar{k}$ such that $k_i=2\La_i+1$. Such summands have a simple pole at $\la_i=\La_i$, see the proof of
Lemma~\ref{la-cont}. In
this case $c_{\bar{m}}(\la,p)$ is equal to zero at $\la_i=\La_i$, so the function
$J_{\bar{l},\bar{m}}(z,\la,p)=c_{\bar{m}}(\la,p)I(w,W)(z,\la,p)$
is well defined at $\la_i=\La_i$. 

This proves Theorem~\ref{an.cont.2}.

Moreover, if $B_\La(\bar{l})=B_\La(\bar{m})=B$, then we have
\be
J_{\bar{l},\bar{m}}(z,\La,p)=\sum_{\bar{k}}\,\frac{l!}{n!\,k_1!\dots
k_n!(l-k)!}\int\limits_{t_i\in{\cal C},\atop
i=k+1,\dots,l}(c_{\bar{m}}\res_{\bar{k}}f)(t_{k+1},\dots,t_l,z,\La,p)\,d^{l-k}t,
\ee
where ${\cal C}$ denotes the curve ${\cal C}_u$ with
$u_m=\min\{-\Real\La_m-\epsilon,0\}$, $m=1,\dots,n$,
and the sum is over all $\bar{k}\in\Z^n_{\ge 0}$ such that
$k\le l$, $k_m\le 2\Real\La_m+1$ for $m=1,\dots,n$, and $k_m=2\La_m+1$ for all $m\in B$.
 
Let $j\in B$.
We compute the residues at points $z_j+\la_j,z_j+\la_j-1,\dots,z_j+\la_j-2\La_j$ explicitly using the
following Lemma.

\begin{lemma}
Let $B_\La(\bar{l})=B_\La(\bar{m})=B$ and let $j\in B$, so that $\La_j\in\La^+$,
$\bar{l},\bar{m}\in\Zb_l^n$ and $l_j,m_j>2\La_j$. 
Let $k=2\La_j$, 
\bea
\lefteqn{(c_{\bar{m}}\res_{t_{k+1}=z_j+\la_j-k}\ldots\res_{t_2=z_j+\la_j-1}\res_{t_1=z_j+\la_j}
\Phi_l w_{\bar{l}}W_{\bar{m}})(t_{k+2},\dots,t_l,z,\La)=}
\\&&
(k+1)!\,e^{\mu(k+1)z_j\pi i/p}\,\psi_{m_j,k}\,
\phi_{\La,j,k}(z)\,(c_{\bar{m}^{(j)}}\Phi_{l-k-1}
w_{\bar{l}^{(j)}}W_{\bar{m}^{(j)}})
(t_{k+2},\dots,t_l,z,\La^{(j)}),
\eea
where $\bar{l}^{(j)}=(l_1,\dots,l_{j-1},l_j-k-1,l_{j+1},\dots,l_n)$,
$\bar{m}^{(j)}=(m_1,\dots,m_{j-1},m_j-k-1,m_{j+1},\dots,m_n)$  and
$\bar{\La}^{(j)}=(\La_1,\dots,\La_{j-1},-\La_j-1,\La_{j+1},\dots,\La_n)$. 
\end{lemma}
\begin{proof}
Consider the case $j=1$.
For a function $g(t,z,\la)$, the coefficient in the Laurent expansion with 
respect to variables $t_1,...,t_{2\La_1+1}$ at $t_i=z_1+\la_1-i+1,\,i=1,\dots,2\La_1+1$, computed at 
$\la_1=\La_1$, i.e. the coefficient of 
\be
\prod_{i=1}^{2\La_1+1}\frac{1}{t_i-z_1-\la_1+i-1},
\ee
is called the \emph{main coefficient}. We compute the main coefficients of the functions
$\Phi_l,w_{\bar{l}},W_{\bar{m}}$.

The main coefficient of the function $\Phi_l(t,z,\la)$ is
\bea
\lefteqn{e^{\mu(2\La_1+1)z_j\pi i/p}\,a_1(\La)\,\phi_{\La^{(1)},1,2\La_1}(z)
\prod_{j=2\La_1+2}^l\prod_{i=1}^{2\La_1+1}
\frac{\Gamma((z_1+\La_1-i+1-t_j-1)/p)}{\Gamma((z_1+\La_1-i+1-t_j+1)/p)}\times}
\\&&
\times\prod_{j=2\La_1+2}^l\prod_{i=1}^n
\frac{\Gamma((t_j-z_i+\la_i)/p)}{\Gamma((t_j-z_i-\la_i)/p)}
\prod_{2\La_1+1< i< j\le l}\frac{\Gamma((t_i-t_j-1)/p)}{\Gamma((t_i-t_j+1)/p)}.
\eea
Here $a_1(\La_1)\in\C$ is some explicitly computable constant. Note that
$\phi_{\La^{(1)},1,2\La_1}(z)=\phi_{\La,1,2\La_1}(z)$.
The main coefficient of the function
$\Phi_l(t,z,\la)$ is equal to
\bea
\lefteqn{\prod_{j=2\La_1+2}^l
\frac{t_j-z_1-\La_1}{t_j-z_1+\La_1}\,\frac{t_j-z_1+\La_1^\prime}{t_j-z_1-\La_1^\prime}
\, \frac{\sin((t_j-z_1-\La_1)\pi i/p)}{\sin((t_j-z_1+\La_1)\pi i/p)}
\,\frac{\sin((t_j-z_1+\La_1^\prime)\pi i/p)}{\sin((t_j-z_1-\La_1^\prime)\pi i/p)}\times}
\\&&
\times a_1(\La)\,e^{\mu(2\La_1+1)z_j\pi i/p}\,\phi_{\La,1,2\La_1}(z)\,
\Phi_{l-2\La_1-1}(t_{2\La_1+1},\dots,t_l,z,\La^{(1)}),
\eea
where $\La_1^\prime=-\La_1-1$.

The function $w_{\bar{l}}(t,z,\la)$ has the form 
$w_{\bar{l}}=\Sum_{\sigma\in {\Bbb S}^l}[\ldots]_\sigma^{rat}$.
Consider the sum over $\sigma\in {\Bbb S}^l$ such that $\sigma(i)=i$ for all 
$i=l_1+1,...,l$. We have 
\bea
\lefteqn{
w_{\bar{l}}(t,z,\la)=
a_2(\La_1)\times}
\\&& \times
\sum_{\sigma \in {\Bbb S}^l}\left[\prod_{1\le i<j\le l_1}\frac{t_i-t_j}{t_i-t_j+1}
\prod_{i=1}^l\frac{1}{t_i-z_1-\la_1}
\prod_{i=l_1+1}^l(t_j-z_1+\la_1)
w_{(l_2,\dots,l_n)}(t_{l_1+1},\dots,t_l,z,\la)
\right] _\sigma^{rat}.
\eea
Here, again, $a_2(\La_1)\in \C$ is an easily computable constant.

Let $\sigma\in {\Bbb S}^l$ be a permutation and there exists $i\le 2\La_1+1$
such that $\sigma(i)>l_1$. Then the term corresponding to $\sigma$ does 
not contribute to the main coefficient. 
Hence, this coefficient is
\bea
\lefteqn{\tilde{a_2}(\La_1)
\sum_{\sigma \in {\Bbb S}^{l-2\La_1-1}}\left[
\prod_{j=2\La_1+2}^{l_1}\prod_{i=1}^{2\La+1}
\frac{z_1+\La_1-i+1-t_j-1}{z_1+\La_1-i+1-t_j+1}
\prod_{i=2\La_1+2}^l\frac{1}{t_i-z_1-\la_1}
\right.\times
}
\\&&
\times\left.
\prod_{i=l_1+1}^l(t_j-z_1+\la_1)
\prod_{2\La_1+2\le i<j\le l_1}\frac{t_i-t_j}{t_i-t_j+1}
w_{(l_2,\dots,l_n)}(t_{l_1+1},\dots,t_l,z,\La_1,\la_2,\dots,\la_n)
\right] _\sigma^{rat},
\eea
where $\tilde{a_2}(\La_1)\in\C$ is another easily computable constant
and the group ${\Bbb
S}^{l-2\La_1-1}$ permutes the variables 
$t_{2\La_1+1},\dots,t_l$. Simplifying, we get the main coefficient of
the function
$w_{\bar{l}}(t,z,\la)$ is given by
\be
\tilde{a_2}(\La_1)
\prod_{j=2\La+1}^l
\frac{t_j-z_1+\La_1}{t_j-z_1-\La_1}\,\frac{t_j-z_1-\La_1^\prime}{t_j-z_1+\La_1^\prime}
\,w_{(l_1^\prime,l_2\dots,l_n)}.
\ee

The main coefficient of the function  $W_{\bar{m}}(t,z,\la)$ is computed similarly and is equal to 
\be
a_3(\La_1)
\prod_{j=2\La_1+1}^l
\frac{\sin((t_j-z_1+\La_1)\pi i/p)}{\sin((t_j-z_1-\La_1)\pi i/p)}
\,
\frac{\sin((t_j-z_1-\La_1^\prime)\pi i/p)}{\sin((t_j-z_1+\La_1^\prime)\pi i/p)}
\,
W_{(m_1^\prime,m_2\dots,m_n)},
\ee
where $a_3(\La_1)\in\C$ is some easily computable comstant.

Multiplying the above main coefficients we get the statement of the Lemma. 
We have
\be
a_1(\La_1)\tilde{a}_2(\La_1)a_3(\La_1)c_{\bar{m}}(\La)=(k+1)!\,
c_{\bar{m}}^{(1)}(\La^{(1)})\,\psi_{2\La_1}.   
\ee

The cases $j=2,\dots,n$  are proved similarly.
\end{proof}

We have
\bean\label{new}
\lefteqn{
J_{\bar{l},\bar{m}}(z,\La,p)=\sum_{\bar{k}}\,\frac{l!}{n!\,k_1!\dots
k_n!(l-k)!}\times}
\\&&
\times\int\limits_{t_i\in{\cal C},\atop
i=k+1+l^\prime(B),\dots,l}
{\frak C}_{\La,\bar{m}}(c_{\bar{m}^\prime(B)}\res_{\bar{k}}
\Phi_{l^\prime(B)}w_{\bar{l}^\prime(B)}W_{\bar{m}^\prime(B)})
(t_{k+1+l^\prime(B)},\dots,t_l,z,\La,p)\,d^{l-k-l^\prime(B)}t,\notag
\eean
where ${\cal C}$ denotes the curve ${\cal C}_u$ with
$u_m=\min\{-\Real\La_m-\epsilon,0\}$, $m=1,\dots,n$,
and the sum is over all $\bar{k}\in\Z^n_{\ge 0}$ such that
$k\le l-l^\prime(B)$, $k_m\le 2\Real\La_m+1$ for $m=1,\dots,n$, and $k_m=0$ for all $m\in B$.

Let ${\cal C}^\prime$ be the curve ${\cal C}_u$ with
$u_m=\min\{-\Real\La_m^\prime-\epsilon,0\}$, $m=1,\dots,n$. Then the integrand of each
integral in sum \Ref{new} does not have poles of the type $t_j=a$, where $a\in\C$ is
located between ${\cal C}^\prime$ and ${\cal C}$, see Lemma~\ref{res}. Thus,
we can move the contour of integration in each summand to 
$\{t\in\C^{l-k-l^\prime(B)}\,|\,t_i\in{\cal C}^\prime,\,i=k+1+l^\prime(B),\dots,l\}$.

We proved Theorem~\ref{functor} for our choice of $z,p$. Theorem~\ref{functor} holds for all $z,p$,
since both left and right hand sides of the formula of Theorem~\ref{functor} are meromorphic functions of
parameters $z,p$.
$\;\Box$

\subsection{Proof of Theorem~\ref{det}.}

Consider a matrix
\be
I^l(z,\la,p)=\{\prod_{i\in B_\La(\bar{m})}(\la_i-\La_i)I_{\bar{l},\bar{m}}(z,\la,p)\}
_{\bar{l},\bar{m}\in\Zb^n_l}.
\ee

The matrix $I^l(z,\la,p)$ is obtained from the matrix $J^l(z,\la)=J^l(z,\la,p)$ by
multiplication the
$\bar{m}$-th row by $(\prod\limits_{i\in B_\La(\bar{m})}(\la_i-\La_i))/c_{\bar{m}}(\la)$ 
for all $\bar{m}\in\Zb^n_l$. For generic 
$p$, these factors are well defined and not equal to zero. 
All entries of the matrix $I^l(z,\la,p)$ are well defined and
the matrix $I^l(z,\la)$ is upper block triangular at $\la=\La$. 

It follows from formula 5.14 in \cite{TV1} that the matrix $I^l(z,\la)$ is non-degenerate for
all
$z,\la,p$ satisfying conditions \Ref{step}-\Ref{weights3}. 

Hence, the matrix \be I^l_{\rm adm}(z,\la,p)=\{\prod_{i\in
B_\La(\bar{m})}(\la_i-\La_i)I_{\bar{l},\bar{m}}(z,\la,p)\}, \ee where $\bar{l},\bar{m}\in\Zb^n_l$ run
through the set of all $\La$-admissible indices, is also non-degenerate for all $z,\la,p$ satisfying
conditions \Ref{step}-\Ref{weights3}. 

For all $p,\La$ satisfying conditions \Ref{step}-\Ref{weights3} and $\La$-admissible
indices $\bar{m}$, the factors
$(\prod\limits_{i\in B_\La(\bar{m})}(\la_i-\La_i))/c_{\bar{m}}(\la)$ are 
well defined and non-zero. Hence, the matrix $J^l_{\rm adm}(z,\la,p)$ is also non-degenerate for
all $z,\la,p$ satisfying conditions \Ref{step}-\Ref{weights3}.

The matrix $J^l(z,\la)$ is upper block triangular at $\la=\La$. 
Now, formula \Ref{det-prod} follows from the Inclusion-Exclusion Principle and
Theorem~\ref{functor}.

The Theorem is proved.$\;\Box$

\subsection{Proof of Theorem~\ref{rational $R$-matrix}}
   
The rational $R$-matrix  $R_{\la_1\la_2}(x)\in\End (V_{\la_1}\T V_{\la_2})$ is a meromorphic function of
$\la_1,\la_2$, defined
explicitly in Section~\ref{rat $R$-matrix} for $\la_1,\la_2\not \in \La^+$.

Consider the qKZ equation with $n=2,\,\mu=2,\,p=-2(|\la_1|+|\la_2|)-1$. Consider the  functions
$\{\Psi_{\bar{m}}(z,\la)\,|\,\bar{m}\in\Zb^2_l\}$ given by
\be
\Psi_{\bar{m}}(z,\la)=\, \sum_{l_1+l_2=l} \,
J_{\bar{l},\bar{m}}(z,\la) \, f^{l_1}v_1\T f^{l_2}v_2,
\ee
where the function $J_{\bar{l},\bar{m}}$ is defined in \Ref{c int}.
   
The functions $\{\Psi_{\bar{m}}(z,\la)\,|\,\bar{m}\in\Zb^2_l\}$ are solutions of the
qKZ equation. Moreover, for generic $\la$, the values of these functions
span $V_{\la_1}\T V_{\la_2}$ over $\C$ by Theorem 5.14 in \cite{TV1}.
By Theorem~\ref{an.cont.2}, these functions are well defined for all
$\la_1,\la_2$. Moreover, they span $V_{\la_1}\T V_{\la_2}$,
for all  $\la_1,\la_2$ and generic $z$, see  theorem 5.14 in
\cite{TV1} and Section~\ref{determinant}.
The first statement of Theorem~\ref{rational $R$-matrix} follows, now, from the equality of meromorphic functions:
\be
\Psi_{\bar{m}}(z_1+p,z_2,\la_1,\la_2)=
e^{-\mu h_1}R_{\la_1\la_2}(z_1-z_2)\Psi_{\bar{m}}(z_1,z_2,\la_1,\la_2)
\ee
for all $\bar{m}\in\Zb^2_l$.

For all complex  $\la_1,\la_2$, the $R$-matrix $R_{\la_1\la_2}(x)$
commutes with the action of $sl_2$. 

If $\la_1,\la_2\not\in\La^+$, then $S_{\la_1}\T V_{\la_2}+V_{\la_1}\T S_{\la_2}$ is the zero submodule
and the second statement of the Theorem is trivial.
Otherwise, let $\la\in\{\la_1,\la_2\}$ be the minimal dominant weight. Namely, if
$\la_1\in\La^+$ and
$\la_2\not\in\La^+$, set $\la=\la_1$. If $\la_2\in\La^+$ and
$\la_1\not\in\La^+$, set $\la=\la_2$. If $\la_1\in\La^+$ and
$\la_2\in\La^+$, set $\la=\min\{\la_1,\la_2\}$.
Then the submodule $S_{\la_1}\T V_{\la_2}+V_{\la_1}\T S_{\la_2}$ is generated by all singular vectors in
$V_{\la_1}\T V_{\la_2}$ of weights less than $\la_1+\la_2-2\la$.
Hence, the $R$-matrix $R_{\la_1\la_2}(x)$ preserves
$S_{\la_1}\T V_{\la_2}+V_{\la_1}\T S_{\la_2}$. 

The third statement of the Theorem follows from the first two statements.
$\;\Box$

\end{document}